\documentclass[10pt,reqno]{extarticle}
\usepackage{amssymb}
\usepackage{amsthm}
\usepackage{amsmath}
\usepackage{colortbl}
\usepackage{float}

\usepackage{pdflscape}

\usepackage[normalem]{ulem}
\usepackage[margin=1in]{geometry}
\usepackage{graphicx}
\usepackage{tikz-qtree}
\usepackage{multicol}
\usepackage{rotating}
\newtheorem{theorem}{Theorem}[section]

\newtheorem{example}[theorem]{Example}

\begin{document}
\begin{center} \Large {\bf Symmetry group classification and optimal reduction of a class of damped Timoshenko beam system
with non-linear rotational moment}
\end{center}
\begin{center} {S. M. Al-Omari$^*$ , F. D. Zaman$^*$, A. Y. Al-Dweik$^*$ , Ryad A Ghanam$^{**}$}
\end{center}
\begin{center}
{$^{*}$Department of Mathematics \& Statistics, King Fahd
University of Petroleum and Minerals,\\ Dhahran 31261, Saudi
Arabia} \\
{$^{**}$ Virginia Commonwealth University in Qatar PO Box 8095,
Doha, Qatar}
\end{center}
\begin{center}
shalomari@kfupm.edu.sa, fzaman@kfupm.edu.sa, aydweik@kfupm.edu.sa
and raghanam@vcu.edu
\end{center}
\begin{abstract}
We consider a non-linear Timoshenko system of partial differential
equations (PDEs) with frictional damping term in rotation angle.
The nonlinearity is due to the arbitrary dependence on the
rotation moment. A Lie symmetry group classification of the
arbitrary function of rotation moment is presented. Optimal system
of one-dimensional subalgebras of the non-linear damped Timoshenko
system are derived for all the non-linear cases. All possible
invariant variables  of the optimal systems for the three
non-linear cases are presented. The corresponding reduced systems
of ordinary differential equations (ODEs) are also provided.
\end{abstract}
\section{Introduction}
The classification of group invariant solutions of differential
equations by means of the optimal systems is one of the main
applications of Lie group analysis to differential equations. The
method was first introduced by Ovsiannikov \cite{fuzel1}. The main
idea behind the method is discussed in his papers
\cite{fuzel2,fuzel3} and also by Chupakin \cite{fuzel4}  and
Ibragimov et al \cite{fuzel5} and Olver \cite{olver}. We can
always construct a family of group invariant solutions
corresponding to a subgroup of a symmetry group admitted by a
given differential equation.  Since there are an infinite number
of such subgroups, it is not possible to list all the group
invariant solutions.  An effective and systematic way of
classifying these group invariant solutions is to obtain optimal
systems of subalgebras of the symmetry Lie algebra. This leads to
non-similar invariant solutions under symmetry transformations.

Timoshenko \cite{Timoshenko} proposed a beam theory which adds the
effect of shear as well as the effect of rotation to the
Euler-Bernoulli beam. The Timoshenko model is a major improvement
for non-slender beams and for high-frequency response where shear
or rotary effects are not negligible \cite{han}. Rivera et. al.
\cite{rivera} studied the global stability for the following
damped Timoshenko beam system with non-linear rotation moment
\begin{equation} \label{1}
\begin{array}{r}
\rho_1 \varphi_{tt} - k (\varphi_x + \psi)_x = 0,\\
\rho_2 \psi_{tt}- (\chi(\psi_x))_x  +k (\varphi_x + \psi) + d \psi_t= 0.\\
\end{array}
\end{equation}
where the functions $\varphi$, $\psi$ depending on $(t,x)\in
(0,\infty)\times (0,L)$ model the transverse displacement of a
beam with $(0,L)\in \Re $, and rotation angle of a filament,
respectively. The constants $\rho_1$, $\rho_2$, $d$ and $k$ are
positive and $\chi$ is a function of $\psi_x$ assumed to satisfy
\begin{equation} \label{2}
\chi_{\psi_{x}}(0)=b,
\end{equation}
with positive constant $b$. However, the algebraic properties of
the Lie algebra admitted by the (\ref{1}) have not been studied so
far. In this paper we perform a Lie symmetry analysis of
non-linear damped Timoshenko beam system (\ref{1}). In sections
two, the complete Lie group classification is presented using
Janet basis. In section three, optimal system of one-dimensional
subalgebras of the non-linear damped Timoshenko system are derived
for all the non-linear cases. In section four, all possible
invariant variables  of the optimal systems for the three
non-linear cases of $\chi(\psi_x)$ are presented. Moreover, the
corresponding reduced systems of ODEs are also provided. As an
illustration, some invariant solutions are given explicitly in
three examples by solving the reduced systems of ODEs.
\section{Complete Lie group classification}
Consider the system of two PDEs with two independent variables
$(t,x)$ and two dependent variables $(\varphi,\psi)$ and the
function $\chi=\chi(\psi_x)$ given by the system (\ref{1}).\\
Following is a brief summary of Lie symmetry
\cite{Ibragimov,Bluman}.

 Consider the following symmetry transformation group acting
  on system of PDEs (\ref{1}).
 \begin{equation} \label{5}
 \begin{array}{ll}
   \tilde{t}=t + \epsilon \xi^1 (t,x,\varphi,\psi)+ O (\epsilon^2),  & \tilde{x}= x + \epsilon \xi^2 (t,x,\varphi,\psi)+ O (\epsilon^2), \\
   \tilde{\varphi}=\varphi + \epsilon \eta^1 (t,x,\varphi,\psi)+ O (\epsilon^2), & \tilde{\psi}=\psi + \epsilon \eta^2 (t,x,\varphi,\psi)+ O(\epsilon^2), \\
 \end{array}
 \end{equation}
where $\epsilon$ is the group parameter and $\xi^1,\xi^2$ and
$\eta^1,\eta^2$ are the infinitesimals of transformations for the
independent and dependent variables, respectively. The associated
 Lie point symmetry generator (vector field) of the system (\ref{1})
 is of the form
\begin{equation} \label{6}
X=\xi^{1}(t,x,\varphi,\psi)\frac{\partial}{\partial
t}+\xi^{2}(t,x,\varphi,\psi)\frac{\partial}{\partial x}+\eta^1
(t,x,\varphi,\psi)\frac{\partial}{\partial \varphi}+\eta^2
(t,x,\varphi,\psi)\frac{\partial}{\partial \psi}.
\end{equation}
The second prolongation of the generator is given by
\begin{equation} \label{7}
X^{[2]}= X + \eta_{i_1}^{\mu} (x^j,u^j,\partial
u^j)\frac{\partial}{\partial u_{{i_1}}^{\mu}} + \eta_{i_1
i_2}^{\mu} (x^j,u^j,\partial u^j,\partial^2 u^j)
\frac{\partial}{\partial u_{{i_1} {i_2}}^{\mu} },
\end{equation}
such that $j=1,2$, and $(x^1,x^2)=(t,x)$, and
$(u^1,u^2)=(\varphi,\psi)$, and $(\partial u^1,\partial
u^2)=(\partial \varphi, \partial \psi)$, and so on... , where
\begin{equation} \label{8}
\eta_{i_1}^{\mu}=D_{i_1} \eta^{\mu}- \sum_{j=1}^2 (D_{i_1} \xi^j)
u_{j}^{\mu}, \text{ }\text{ } \mu=1,2 \text{ } , \nonumber
\end{equation}
\begin{equation} \label{9}
\eta_{i_1 i_2}^{\mu}=D_{i_2} \eta_{i_1}^{\mu}- \sum_{j=1}^2
(D_{i_2} \xi^j) u_{i_1 j}^{\mu},  \nonumber
\end{equation}
\begin{equation} \label{10}
D_i=\frac{\partial}{\partial x_i}+ u_{i}^{\mu}
\frac{\partial}{\partial u^{\mu}}+ u_{i j}^{\mu}
\frac{\partial}{\partial u_{j}^{\mu}}+ u_{i i_1 i_2}^{\mu}
\frac{\partial}{\partial u_{i_1 i_2}^{\mu}}+... ,
\end{equation}
where $D_i$ is the total derivative operator.\\

Using the invariance condition of the system of PDEs (\ref{1})
\begin{equation} \label{11}
\begin{array}{l}
X^{[2]} (\rho_1 \varphi_{tt} - k (\varphi_x + \psi)_x)\vert_{(\ref{1})}=0,\\
X^{[2]} (\rho_2 \psi_{tt}- \chi_x (\psi_x) +k (\varphi_x + \psi) +d  \psi_t) \vert_{(\ref{1})}=0,\\
\end{array}
\end{equation}
and comparing coefficients of the various derivatives of the
dependent variables $\phi$ and $\psi$ yields an over-determined
linear PDE system. Carrying out the Janet basis of this
over-determined system in the degree reverse lexicographical
ordering as $\psi>\phi>x>t$ and $ \eta_2>\eta_1>\xi_2>\xi_1$   by
using the command "JanetBasis" involved in the Maple package
"Janet" \cite{Blinkov2003}, leads to two cases. These two cases
arise from the command "Denominators" which returns the functions
by which the Janet basis algorithm had to divide. These two cases
are given as follows:
\subsection{ $\chi(\psi_x)=b \psi_x + \gamma $}
In this case $\chi(\psi_x)$ is a linear function satisfying the
condition (\ref{2}). The Janet basis of  of the over-determined
system is
\begin{equation}\label{12}
    \begin{array}{l}
      [\eta^1_\psi,\xi^2_\psi,\xi^1_\psi,\eta^2_\phi,\eta^1_\phi-\eta^2_\psi,\xi^2_\phi,
      \xi^1_\phi,\xi^2_x,\xi^1_x,\xi^2_t,\xi^1_t,\eta^2_{\psi \psi},\eta^2_{x \psi},
      \eta^2_{t \psi},\eta^1_{t \psi},\eta^2_{x \phi}, \\ \eta^2_{t \phi},\eta^1_{t \phi},
      \frac {k}{\rho2}\eta^1_x+\frac {k}{\rho_2}\eta^2+\frac{d}{\rho_2}\eta^2_t-\frac{k}{\rho_2}\eta^2_{\psi} \psi-\frac {b}{\rho_2}\eta^2_{x,x}+\eta^2_{t,t},
      -\frac{k}{\rho_1}\eta^1_{x,x}+\eta^1_{t,t}-\frac{k}{\rho_1}\eta^2_x].\\
    \end{array}
\end{equation}
The solution of this system of determining equations is
\begin{equation}\label{}
 \begin{array}{cccc}
   \xi^1 = c_1, & \xi^2=c_2, & \eta^1=  c_3 \varphi + f(t,x), & \eta^2= c_3 \psi + g(t,x), \\
 \end{array}
\end{equation}
where $f(t,x)$ and $g(t,x)$ satisfy the following system of PDEs
\begin{equation}\label{}
    \begin{array}{r}
      k f_x + k g + d g_t - b g_{xx} + \rho_2 g_{tt}=0, \\
      \rho_1 f_{tt} - k f_{xx} - k g_x =0.\\
    \end{array}
\end{equation}
The corresponding Lie point symmetry generators admitted by the
system (\ref{1}) are given as
\begin{equation} \label{13}
\begin{array}{cccc}
  X_1=  \frac{\partial}{\partial t}, & X_2=  \frac{\partial}{\partial x}, & X_3= \varphi \frac{\partial}{\partial \varphi}+\psi \frac{\partial}{\partial \psi}, &  X_{\infty}=f(t,x) \frac{\partial}{\partial \varphi}+ g(t,x) \frac{\partial}{\partial \psi}, \\
\end{array}
\end{equation}

\subsection{$\chi_{\psi_x \psi_x} \neq 0 $}
The Janet basis of the over-determined system is
\begin{equation}\label{15}
\begin{array}{l}
[\eta^2_{{\psi_{{}}}},\eta^1_{{\psi_{{}}}},\xi^2_{{\psi_{{}}}},\xi^1_{{\psi_{{}}}},\eta^2_{{\phi_{{}}}},\eta^1_{{\phi_{{}}}},\xi^2_{{\phi_{{}}}},\xi^1_{{\phi_{{}}}},\eta^2_{{x}},\xi^2_{{x}},\xi^1_{{x}},\xi^2_{{t}},\xi^1_{{t}},\\
\eta^1_{{x,\psi_{{}}}},\eta^2_{{t,\psi_{{}}}},\eta^1_{{t,\psi_{{}}}},\eta^1_{{x,\phi_{{}}}},\eta^2_{{t,\phi_{{}}}},\eta^1_{{t,\phi_{{}}}},\eta^1_{{x,x}},\eta^2_{{t,x}},{\frac{k}{ \rho_{{2}}} \eta^1_{{x}} }+{\frac{k}{\rho_{{2}}}\eta^2_{{}}}+{\frac {d}{\rho_{{2}}}\eta^2_{{t}}}+\eta^2_{{t,t}},\eta^1_{{t,t}},\\
\eta^1_{{t,x,\psi_{{}}}},\eta^1_{{t,x,\phi_{{}}}},\eta^1_{{t,x,x}}]\\
\end{array}
\end{equation}
The solution of the system (\ref{15}) is
\begin{equation}\label{16}
 \begin{array}{cccc}
   \xi^1 = c_1, & \xi^2=c_2, & \eta^1=  c_3+ c_4~t+ c_5 x+c_6~tx, & \eta^2= F(t) \\
 \end{array}
\end{equation}
where $F(t)$ satisfies the following ODE
\begin{equation} \label{17}
\rho_2 F''(t) + d F'(t) + k F(t) = - c_6 k t - c_5 k.
\end{equation}
The characteristic equation of left hand side of the equation
(\ref{17}) gives rise to the following three sub-cases
\subsubsection{$d^2-4 k\rho_2 = 0$}
The Lie point symmetry generators admitted by the system (\ref{1})
are given by
\begin{equation}\label{18}
\begin{array}{llll}
 X_1=  \frac{\partial}{\partial t} & X_2=  \frac{\partial}{\partial x},
 &X_3=  \frac{\partial}{\partial \varphi},\\
 X_4=t \frac{\partial}{\partial \varphi},&
 X_5= x \frac{\partial}{\partial \varphi}- \frac{\partial}{\partial \psi},
 & X_6=tx \frac{\partial}{\partial \varphi} + (2\sqrt{\frac{\rho_2}{k}}-t) \frac{\partial}{\partial
 \psi},\\
 X_7=e^{-\sqrt{\frac{k}{\rho_2}}t}\frac{\partial}{\partial \psi},
 &X_8= t e^{-\sqrt{\frac{k}{\rho_2}}t}  \frac{\partial}{\partial \psi}.\\
\end{array}
\end{equation}
In order to obtain the group transformations which are generated
by the resulting infinitesimal symmetry generators (\ref{18}), we
need to solve the following system of first order ODEs
\begin{equation} \label{19}
\begin{array}{lll}
  \frac{d \tilde{x}^j (\epsilon)}{d \epsilon} = \xi^j (\tilde{t}(\epsilon),\tilde{x}(\epsilon),\tilde{\varphi}(\epsilon),\tilde{\psi}(\epsilon)), & \tilde{x}^j  (0)=x^j,\\
  \frac{d \tilde{u}^j (\epsilon)}{d \epsilon} = \eta^j (\tilde{t}(\epsilon),\tilde{x}(\epsilon),\tilde{\varphi}(\epsilon),\tilde{\psi}(\epsilon)), & \tilde{u}^j (0)=u^j, &j=1,2.\\
\end{array}
\end{equation}
The one parameter group $G_i(\epsilon)=e^{\epsilon X_i}$ generated
by $X_i$ for $i=1,...,8$, are as follows:
\begin{equation}\label{20}
\begin{array}{ll}
 G_1(\epsilon):(t,x,\varphi,\psi)\mapsto (t+\epsilon, x, \varphi , \psi),   & G_5(\epsilon):(t,x,\varphi,\psi)\mapsto (t, x, \varphi+\epsilon x , \psi-\epsilon),  \\
 G_2(\epsilon):(t,x,\varphi,\psi)\mapsto (t, x+\epsilon, \varphi , \psi),   & G_6(\epsilon):(t,x,\varphi,\psi)\mapsto (t, x, \varphi+\epsilon t x , \psi+ \epsilon (2 \sqrt{\frac{\rho_2}{k}}-t)), \\
 G_3(\epsilon):(t,x,\varphi,\psi)\mapsto (t, x, \varphi+\epsilon , \psi),   & G_7(\epsilon):(t,x,\varphi,\psi)\mapsto (t, x, \varphi , \psi+\epsilon e^{-\sqrt{\frac{k}{\rho_2}}t}), \\
 G_4(\epsilon):(t,x,\varphi,\psi)\mapsto (t, x, \varphi+\epsilon t , \psi), & G_8(\epsilon):(t,x,\varphi,\psi)\mapsto (t, x, \varphi , \psi+\epsilon t e^{-\sqrt{\frac{k}{\rho_2}}t}).\\
\end{array}
\end{equation}
\begin{theorem}
If $\varphi=f(t,x)$ and $\psi=g(t,x)$ is a solution of the
Timoshenko system (\ref{1}) with $d^2 - 4 k \rho_2 =0$, then so is
 \begin{equation}\label{30-3}
    \begin{array}{l}
      \varphi=f(t+\epsilon_1, x + \epsilon_2) + \epsilon_3 + \epsilon_4 ( t+\epsilon_1) +\epsilon_5( x+\epsilon_2)
      + \epsilon_6( t+\epsilon_1)(x+\epsilon_2), \\
      \psi =g(t+\epsilon_1, x + \epsilon_2) + 2 \sqrt{\frac {\rho_2}{k}} \epsilon_6 - \epsilon_6
      (t+\epsilon_1) -\epsilon_5 + \epsilon_7 e^{-\sqrt {\frac {k}{\rho_2}}( t+\epsilon_1) }+ \epsilon_8 ( t+\epsilon_1) e^{-\sqrt {\frac {k}{\rho_2}}( t+\epsilon_1) },\\
    \end{array}
\end{equation}
where $\{\epsilon_i\}_{i=1}^8$ are arbitrary real numbers.
\begin{proof}
The eight parameters group
$$G(\epsilon_1,\epsilon_2,\epsilon_3,\epsilon_4,\epsilon_5,\epsilon_6,\epsilon_7,\epsilon_8)=G_8(\epsilon_8)\circ
G_7(\epsilon_7) \circ G_6(\epsilon_6)\circ G_5(\epsilon_5)\circ
G_4(\epsilon_4)\circ G_3(\epsilon_3)\circ G_2(\epsilon_2)\circ
G_1(\epsilon_1)$$ generated by $X_i$ for $i=1,...,8$, can be given
by the composition of the transformations (\ref{20}) as follows:
\begin{equation}\label{30-2}
    \begin{array}{ll}
    G: \left(t,x,\varphi , \psi\right) \longmapsto &(t+\epsilon_1 ,  x+\epsilon_ 2, \varphi+\epsilon_3+
\epsilon_4 (t+\epsilon_1) +\epsilon_5(x+\epsilon_2)+\epsilon_6
(t+\epsilon_1)(x+\epsilon_2 ),\\& \psi+ 2 \sqrt {\frac
{\rho_2}{k}} \epsilon_6 - \epsilon_6(t+\epsilon_1) -\epsilon_5+
e^{-\sqrt {\frac {k}{\rho_2}}( t+\epsilon_1) }(\epsilon_7
+ \epsilon_8 ( t+\epsilon_1) )), \\
    \end{array}
\end{equation}
and this completes the proof.
\end{proof}
\end{theorem}
\subsubsection{$ d^2 - 4 k \rho_2 = \lambda^2$ such that $\lambda > 0$}
The Lie point symmetry generators admitted by the system (\ref{1})
are given by
\begin{equation}\label{22}
\begin{array}{llll}
 X_1=  \frac{\partial}{\partial t}, & X_2=  \frac{\partial}{\partial x},
 &X_3=  \frac{\partial}{\partial \varphi},\\
 X_4=t \frac{\partial}{\partial \varphi}, &
 X_5= x \frac{\partial}{\partial \varphi}- \frac{\partial}{\partial \psi},
 & X_6=tx \frac{\partial}{\partial \varphi} + (\frac{\sqrt{\lambda^2+4 k \rho_2}}{k}-t) \frac{\partial}{\partial
 \psi},\\
 X_7=e^{-\frac{{\sqrt{\lambda^2+4 k \rho_2}}}{2 \rho_2}t} \cosh (\frac{\lambda t}{2\rho_2}) \frac{\partial}{\partial
 \psi},
 &X_8= e^{-\frac{{\sqrt{\lambda^2+4 k \rho_2}}}{2\rho_2}t}\sinh (\frac{\lambda t}{2\rho_2})  \frac{\partial}{\partial \psi}.\\
\end{array}
\end{equation}
The one parameter group $G_i(\epsilon)=e^{\epsilon X_i}$ generated
by $X_i$ for $i=1,...,8$ are as follows:\\
$G_i(\epsilon), i=1,...,5$ are the same as in equation (\ref{20}),
and $G_i(\epsilon), i=6,7,8$ are given by
\begin{equation}  \label{23}
\begin{array}{l}
   G_6(\epsilon):(t,x,\varphi,\psi)\mapsto (t, x, \varphi+\epsilon t x , \psi+ \epsilon (\frac{{\sqrt{\lambda^2+4 k \rho_2}}}{k}-t), \\
   G_7(\epsilon):(t,x,\varphi,\psi)\mapsto (t, x, \varphi , \psi+\epsilon e^{-\frac{{\sqrt{\lambda^2+4 k \rho_2}}}{\rho_2}t} \cosh (\frac{\lambda t}{2 \rho_2})), \\
   G_8(\epsilon):(t,x,\varphi,\psi)\mapsto (t, x, \varphi , \psi+\epsilon e^{-\frac{{\sqrt{\lambda^2+4 k \rho_2}}}{\rho_2}t} \sinh (\frac{\lambda t}{2 \rho_2})). \\
\end{array}
\end{equation}
\begin{theorem}
If $\varphi=f(t,x)$ and $\psi=g(t,x)$ is a solution of the
Timoshenko system (\ref{1}) with $d^2 - 4 k \rho_2 = \lambda^2$,
then so is
  \begin{equation}\label{30-6}
    \begin{array}{ll}
      \varphi=&F(t+\epsilon_1, x + \epsilon_2) + \epsilon_3+\epsilon_4(t+\epsilon_1)+\epsilon_5 (x+\epsilon c_2)+\epsilon_6
      (t+\epsilon_1)(x+\epsilon_2), \\
      \psi =&G(t+\epsilon_1, x + \epsilon_2)-\epsilon_6(t+\epsilon_1)-\epsilon_5+\frac{d}{k}\epsilon_6  +  e^{-\frac {(d t+\epsilon_1)}{2\rho_2 }}(\epsilon_7 \cosh
     (\frac {\lambda(t+\epsilon_1)}{2\rho_2})+  \epsilon_8 \sinh(\frac {\lambda(t+\epsilon_1)}{2\rho_2})),  \\
    \end{array}
\end{equation}
where $\{\epsilon_i\}_{i=1}^8$ are arbitrary real numbers.
\begin{proof}
The eight parameters group
$$G(\epsilon_1,\epsilon_2,\epsilon_3,\epsilon_4,\epsilon_5,\epsilon_6,\epsilon_7,\epsilon_8)=G_8(\epsilon_8)\circ
G_7(\epsilon_7) \circ G_6(\epsilon_6)\circ G_5(\epsilon_5)\circ
G_4(\epsilon_4)\circ G_3(\epsilon_3)\circ G_2(\epsilon_2)\circ
G_1(\epsilon_1)$$ generated by $X_i$ for $i=1,...,8$, can be given
by the composition of the transformations (\ref{23}) as follows:
\begin{equation}\label{30-5}
    \begin{array}{lll}
     G:(t,x,\varphi , \psi) \longmapsto &(t+\epsilon_1 , x+\epsilon_2 ,\varphi+\epsilon_3+
     \epsilon_4(t+\epsilon_1)+\epsilon_5 (x+\epsilon_2)+ \epsilon_6(t+\epsilon_1)(x+\epsilon_2) ,
       \\ & \psi -\epsilon_6(t+\epsilon_1)-\epsilon_5+\frac{d}{k}\epsilon_6 +e^{-\frac {(d t+\epsilon_1)}{2\rho_2 }}( \epsilon_7 \cosh
     (\frac {\lambda(t+\epsilon_1)}{2\rho_2})+\epsilon_8 \sinh(\frac {2\lambda(t+\epsilon_1)}{2\rho_2}))). \\
    \end{array}
\end{equation}
and this completes the proof.
\end{proof}
\end{theorem}
\subsubsection{$ d^2 - 4 k \rho_2 = - \mu^2 $, such that $\mu > 0$}
The Lie point symmetry generators admitted by the system (\ref{1})
are given by
\begin{equation}\label{25}
\begin{array}{llll}
 X_1=  \frac{\partial}{\partial t}, & X_2=  \frac{\partial}{\partial x},
 &X_3=  \frac{\partial}{\partial \varphi},\\
 X_4=t \frac{\partial}{\partial \varphi},&
 X_5= x \frac{\partial}{\partial \varphi}- \frac{\partial}{\partial \psi},
 & X_6=tx \frac{\partial}{\partial \varphi} + (\frac{\sqrt{4 k \rho_2-\mu^2}}{k}-t) \frac{\partial}{\partial
 \psi},\\
 X_7= e^{-\frac{{\sqrt{4 k \rho_2-\mu^2}}}{2 \rho_2}t} \cos (\frac{\mu t}{2 \rho_2})\frac{\partial}{\partial \psi},
 &X_8=  e^{-\frac{{\sqrt{4 k \rho_2-\mu^2}}}{2 \rho_2}t}
  \sin (\frac{\mu t}{2 \rho_2})
  \frac{\partial}{\partial \psi}.\\
\end{array}
\end{equation}
The one parameter group $G_i(\epsilon)=e^{\epsilon X_i}$ generated
by $X_i$ for $i=1,...,8$ are as follows:\\
$G_i(\epsilon), i=1,...,5$ are the same as in equation (\ref{20}),
and $G_i(\epsilon), i=6,7,8$ are given by
\begin{equation}  \label{26}
\begin{array}{l}
  G_6(\epsilon):(t,x,\varphi,\psi)\mapsto (t, x, \varphi+\epsilon t x , \psi+ \epsilon (\frac{d}{k}-t), \\
  G_7(\epsilon):(t,x,\varphi,\psi)\mapsto (t, x, \varphi , \psi+\epsilon e^{-\frac{d}{\rho_2}t} \cos (\frac{\mu t}{2 \rho_2})), \\
  G_8(\epsilon):(t,x,\varphi,\psi)\mapsto (t, x, \varphi , \psi+\epsilon e^{-\frac{d}{\rho_2}t} \sin (\frac{\mu t}{2 \rho_2})). \\
\end{array}
\end{equation}
\begin{theorem}
If $\varphi=f(t,x)$ and $\psi=g(t,x)$ is a solution of the
Timoshenko system (\ref{1}) with $d^2 - 4 k \rho_2 = - \mu^2$,
then so is
  \begin{equation}\label{30-12}
    \begin{array}{l}
      \varphi=F(t+\epsilon_1, x + \epsilon_2) + \epsilon_3+
     \epsilon_4(t+\epsilon_1)+\epsilon_5(x+\epsilon_2)+\epsilon_6(t+ \epsilon_1)(x+\epsilon_2) \\
      \psi =G(t+\epsilon_1, x + \epsilon_2)-\epsilon_6 (t+\epsilon_1)-\epsilon_5+ \frac{d}{k} \epsilon_6
     +  e^{-\frac {d(t+\epsilon_1)}{2 \rho_2}}( \epsilon_7 \cos(\frac{\mu(t+\epsilon_1)}{2 \rho_2}) +
     \epsilon_8 \sin(\frac {\mu ( t+\epsilon_1)}{2 \rho_2})).  \\
    \end{array}
\end{equation}
where $\{\epsilon_i\}_{i=1}^8$ are arbitrary real numbers.
\begin{proof}
The eight parameter group
$$G(\epsilon_1,\epsilon_2,\epsilon_3,\epsilon_4,\epsilon_5,\epsilon_6,\epsilon_7,\epsilon_8)=G_8(\epsilon_8)\circ
G_7(\epsilon_7) \circ G_6(\epsilon_6)\circ G_5(\epsilon_5)\circ
G_4(\epsilon_4)\circ G_3(\epsilon_3)\circ G_2(\epsilon_2)\circ
G_1(\epsilon_1)$$ generated by $X_i$ for $i=1,...,8$, can be given
by the composition of the transformations (\ref{26}) as follows:
\begin{equation}\label{30-11}
    \begin{array}{lll}
     G:(t,x,\varphi , \psi) \longmapsto &(t+\epsilon_1 , x+\epsilon_2, \phi+\epsilon_3+
     \epsilon_4(t+\epsilon_1)+\epsilon_5(x+\epsilon_2)+\epsilon_6(t+ \epsilon_1)
     (x+\epsilon_2) ,\\& \psi-\epsilon_6 (t+\epsilon_1)-\epsilon_5+ \frac{d}{k} \epsilon_6
      +  e^{-\frac {d(t+\epsilon_1)}{2 \rho_2}}( \epsilon_7 \cos(\frac{\mu(t+\epsilon_1)}{2 \rho_2}) +
     \epsilon_8 \sin(\frac {\mu ( t+\epsilon_1)}{2 \rho_2})). \\
    \end{array}
\end{equation}
and this completes the proof.
\end{proof}
\end{theorem}
\section{Optimal system of one-dimensional subalgebras of the non-linear damped Timoshenko system}
In this section, we give the complete classification of the
one-dimensional optimal system for each of the algebras with basis
(\ref{18}), (\ref{22}) and (\ref{25}). In order to find the
optimal system, one needs to classify the one-dimensional
subalgebras under the action of the adjoint representation. We
follow the algorithm explained by Olver \cite{olver}.\\
First, we calculate the adjoint representation given by
\begin{equation}\label{3-1}
    \begin{array}{l}
      Ad(exp(\epsilon  X_i). X_j)=X_j-\epsilon [X_i , X_j]+\frac{\epsilon^2}{2!}[X_i,[X_i,X_j]]-\frac{\epsilon^3}{ 3 ! }[X_i,[X_i,X_j]]+ ... ,
      \\ \nonumber
    \end{array}
\end{equation}
and then establish the adjoint table. For each case, we classify
the conjugacy classes under the adjoint representation according
to the sign of the killing form.
\subsection {Optimal system for the case $d^2-4 k\rho_2 = 0$}
The non-zero commutators of the Lie algebra $\mathcal{L}^8$ with
basis (\ref{18}) are given by
\begin{eqnarray} \label{21}
\begin{tabular}{llllll}
$ [X_1,X_4]=X_3 $,&$ [X_1,X_6]=X_5 $ ,&$[X_1,X_7]= - \sqrt{\frac{k}{\rho_2}}X_7$, \\
$[X_1,X_8]=X_7 - \sqrt{\frac{k}{\rho_2}}X_8 $,& $[X_2,X_5]=X_3$,&$[X_2,X_6]=X_4$.\\
\end{tabular}
\end{eqnarray}
The Lie algebra $\mathcal {L}^8$ is solvable and the Killing form
is given by $K=2~\hat{a}^2 a^2_1$ where
$\hat{a}=\sqrt{\frac{k}{\rho_2}}$. The adjoint table is given by
\begin{table} [H]  \label{tab:1}
\begin{center}  \caption {}
\begin{tiny}
\begin{tabular}{ccccccccc}
  \hline
 \rowcolor[gray] {0.9} $Ad (e^\epsilon)$ & $X_1$ & $X_2$ & $X_3$ & $X_4$ & $X_5$ & $X_6$ & $X_7$ & $X_8$ \\
  \hline \hline
  $X_1$ & $X_1$ & $X_2$ & $X_3$ & $X_4-\epsilon X_3$ & $X_5$ & $X_6-\epsilon X_5$ & $e^{\epsilon \hat{a}}X_7$ & $e^{\epsilon \hat{a}}X_8-\epsilon e^{\epsilon \hat{a}}X_7 $ \\\hline
  $X_2$ & $X_1$ & $X_2$ & $X_3$ & $X_4$ & $X_5-\epsilon X_3$ & $X_6-\epsilon X_4$ & $X_7$ & $X_8$ \\\hline
  $X_3$ & $X_1$ & $X_2$ & $X_3$ & $X_4$ & $X_5$ & $X_6$ & $X_7$ & $X_8$ \\\hline
  $X_4$ & $X_1+\epsilon X_3$ & $X_2$ & $X_3$ & $X_4$ & $X_5$ & $X_6$ & $X_7$ & $X_8$ \\\hline
  $X_5$ & $X_1$ & $X_2+\epsilon X_3$ & $X_3$ & $X_4$ & $X_5$ & $X_6$ & $X_7$ & $X_8$ \\\hline
  $X_6$ & $X_1+\epsilon X_5$ & $X_2+\epsilon X_4$ & $X_3$ & $X_4$ & $X_5$ & $X_6$ & $X_7$ & $X_8$ \\\hline
  $X_7$ & $X_1-\hat{a} \epsilon X_7$ & $X_2$ & $X_3$ & $X_4$ & $X_5$ & $X_6$ & $X_7$ & $X_8$ \\\hline
  $X_8$ & $X_1+\epsilon X_7-\hat{a} \epsilon X_8$ & $X_2$ & $X_3$ & $X_4$ & $X_5$ & $X_6$ & $X_7$ & $X_8$ \\\hline
\end{tabular}
\end{tiny}
\end{center}
\end{table}
The adjoint group is defined by the matrix
\begin{equation}\label{3-2}
    \begin{array}{c}
      A=Ad(e^{-\epsilon_8 X_8}). Ad(e^{-\epsilon_7 X_7}).
Ad(e^{-\epsilon_6 X_6}).  Ad(e^{-\epsilon_5 X_5}).
Ad(e^{-\epsilon_4 X_4}).  Ad(e^{-\epsilon_3 X_3}).
Ad(e^{-\epsilon_2 X_2}).
Ad(e^{-\epsilon_1 X_1}), \nonumber \\
    \end{array}
\end{equation}
which is given by
\begin{scriptsize}
\begin{equation} \label{3-2-1}
 A= \begin{pmatrix}
  \begin{array}{cccccccc}
  1 & 0 & -\epsilon_4 & 0 & -\epsilon_6 & 0 & \hat{a} \epsilon_7-\epsilon_8 & \hat{a} \epsilon_8 \\
  0 & 1 & -\epsilon_5 & -\epsilon_6 & 0 & 0 & 0 & 0 \\
  0 & 0 & 1 & 0 &0 &0 & 0 & 0 \\
  0 & 0 & \epsilon_1 & 1 & 0 & 0 & 0 & 0 \\
  0 & 0 & \epsilon_2 & 0 & 1 & 0& 0 & 0 \\
  0 & 0 & \epsilon_1 \epsilon_2 & \epsilon_2 & \epsilon_1 & 1 & 0 & 0 \\
  0 & 0 & 0 & 0 & 0 & 0 & e^{-\hat{a} \epsilon_1} & 0 \\
 0 & 0 & 0 & 0 & 0 & 0 &  \epsilon_1 e^{-\hat{a} \epsilon_1} & e^{-\hat{a} \epsilon_1} \\
\end{array}
\end{pmatrix}.
\end{equation}
\end{scriptsize}

\begin{theorem} An optimal system of one-dimensional Lie algebra $\mathcal{L}^8$ with basis (\ref{18}) is provided by the following
generators
\begin{scriptsize}
\begin{equation}    \label{28}
\begin{array}{ll|lll}
  X^1=X_1 + \alpha X_2 + \beta X_6, &  \alpha \in \mathbb{R},\beta \neq 0,                   & X^8= \alpha X_4 + X_5+\beta X_8, &  \alpha \in \mathbb{R}, \beta \neq 0,  \\
  X^2=X_1 + \alpha X_2 + \beta X_4, & \alpha,\beta \in \mathbb{R},                           & X^{9}=\alpha X_4 + \beta X_5 + X_7, & \alpha \in \mathbb{R}, \beta \neq 0, \\
  X^3=X_2 + \alpha X_6 + \beta X_7+ \gamma X_8, & \alpha \neq 0,\beta,\gamma \in \mathbb{R}, & X^{10}=\alpha X_4 + X_5 , & \alpha \in \mathbb{R},   \\
  X^4=X_2 + \alpha X_5 + \beta X_8, & \alpha \in \mathbb{R},\beta \neq 0,                    & X^{11}= X_4 + \alpha X_7 + \beta X_8, & \alpha,\beta \in \mathbb{R},  \\
  X^5=\alpha X_2 + \beta X_5 +  X_7, & \alpha \neq 0 , \beta \in \mathbb{R},                  & X^{12}=  \alpha X_3 + X_8, &  \alpha \in \mathbb{R}, \\
  X^6=X_2 + \alpha X_5 , & \alpha \in \mathbb{R},                                            & X^{13}= \alpha X_3 + X_7, & \alpha \in \mathbb{R},  \\
  X^7=\alpha X_3 + X_6 + \beta X_7 + \gamma X_8, & \alpha,\beta,\gamma \in \mathbb{R},       & X^{14}=X_3.  \\
\end{array}
\end{equation}
\end{scriptsize}

\begin{proof}
Let $X$ and $\tilde{X}$ be two elements in the Lie algebra
$\mathcal{L}^8$ with basis (\ref{18}) given by $X=\sum_{i=1}^{8}
a_i X_i$ and $ \tilde{X}=\sum_{i=1}^{8} \tilde{a}_i X_i$. For
simplicity, we will write $X$ and $\tilde{X}$ as row vectors of
the coefficients on the form $X=(\,a_1 \,\, a_2 \,\, . \,\, . \,\,
. \,\, a_8\,)$ and $\tilde{X}=(\,\tilde{a_1} \,\, \tilde{a_2} \,\,
. \,\, . \,\, . \,\, \tilde{a_8}\,)$. Then in order for $X$ and
$\tilde{X}$ to be in the same conjugacy class, we must have
$\tilde{X}=XA$, where $A$ is given by (\ref{3-2-1}). So, the
theorem is proved by solving the system

\begin{equation}\label{3-4}
 \begin{array}{l}
   \tilde{a_1}= a_1,\\
   \tilde{a_2}= a_2, \\
   \tilde{a_3}= a_3+\epsilon_1 a_4+\epsilon_2 a_5-\epsilon_5 a_2 -\epsilon_4 a_1+\epsilon_1 \epsilon_2 a_6,\\
   \tilde{a_4}= a_4-\epsilon_6 a_2 +\epsilon_2 a_6, \\
   \tilde{a_5}= a_5-\epsilon_6 a_1+\epsilon_1a_6, \\
   \tilde{a_6}= a_6, \\
   \tilde{a_7}= a_1(\hat{a} \epsilon_7-\epsilon_8)+ e^{-\epsilon_1 \hat{a}}(a_7+a_8 \epsilon_1), \\
   \tilde{a_8}= a_1 \hat{a} \epsilon_8  + a_8 e^{-\epsilon_1 \hat{a}}, \\
 \end{array}
\end{equation}

for $\{\epsilon_i\}_{i=1}^8$ in term of $\{a_i\}_{i=1}^8$ in order
to get the simplest values of $\{\tilde{a}_i\}_{i=1}^8$. The
results are presented for different cases in the tree diagram (1)
given in the appendix where it is initiated by the sign of the
Killing form and its leafs are given completely.
The full details for each leaf are given as follows:  \\
\begin{scriptsize}
\uline{\textbf{Case 1}} $a_1 \neq 0$, $a_6 \neq 0$ : Let
$\epsilon_1=\epsilon_5=0$, $\epsilon_2=\frac{a_2 a_5 - a_1
a_4}{a_1 a_6} $, $\epsilon_4=\frac{a_1 a_3 a_6-a_1 a_4 a_5 - a_2
a_5}{a_1^2 a_6}$, $\epsilon_6=\frac{a_5}{a_1}$,
$\epsilon_7=-\frac{\gamma a_7 +a_8}{\gamma^2 a_1^2 }$ and
$\epsilon_8=-\frac{a_8}{\gamma a_1}$ to have
$\tilde{a_3}=\tilde{a_4}=\tilde{a_5} =\tilde{a_7}=\tilde{a_8}=0$.
Then the conjugacy class is $\langle X_1+\alpha X_2+\beta X_6
\rangle$, with $\alpha \in \mathbb{R},\beta \neq 0$.\\
\uline{\textbf{Case 2}} $a_1 \neq 0$, $a_6 = 0$ : Let
$\epsilon_1=\epsilon_2=\epsilon_5=0$,
$\epsilon_4=\frac{a_3}{a_1}$, $\epsilon_6=\frac{a_5}{a_1}$,
$\epsilon_7=-\frac{\gamma a_7+a_8}{\gamma^2 a_1}$,
$\epsilon_8=-\frac{a_8}{\gamma a_1}$ to have
$\tilde{a_3}=\tilde{a_5} =\tilde{a_7}=\tilde{a_8}=0$. So the
conjugacy class is $\langle X_1+\alpha X_2+\beta X_4
\rangle$, with $\alpha,\beta \in \mathbb{R}$.\\
\uline{\textbf{Case 3}} $a_1=0$, $a_2 \neq 0$, $a_6 \neq 0$ : Let
$\epsilon_1=-\frac{a_5}{a_6}$, $\epsilon_2=0$,
$\epsilon_5=\frac{a_3 a_6-a_4 a_5}{a_2 a_6}$, and
$\epsilon_6=\frac{a_4}{a_2}$ to make
$\tilde{a_3}=\tilde{a_4}=\tilde{a_5}=0$. Then the conjugacy class
is $\langle X_2 +\alpha X_6 + \beta X_7+\gamma X_8 \rangle$,
$\alpha \neq 0, \beta,\gamma \in \mathbb{R}$.\\
\uline{\textbf{Case 4}} $a_1=0$, $a_2 \neq 0$, $a_6= 0$, $a_8 \neq
0$ : Let $\epsilon_1=-\frac{a_7}{a_8}$, $\epsilon_2=0$,
$\epsilon_5=\frac{a_3 a_8- a_4 a_7}{a_2 a_8}$, and
$\epsilon_6=\frac{a_4}{a_2}$ to make
$\tilde{a_3}=\tilde{a_4}=\tilde{a_7}=0$. Then the conjugacy class
is of the form $\langle X_2
+\alpha X_5 + \beta X_8 \rangle$, $\alpha \in \mathbb{R},\beta \neq 0$.\\
\uline{\textbf{Case 5}} $a_1=0$, $a_2 \neq 0$, $a_6= 0$, $a_8= 0$,
$a_7 \neq 0$ : Let $\epsilon_1=\frac{\ln\mid a_7 \mid}{\gamma}$,
$\epsilon_2=0$, $\epsilon_5=\frac{\frac{a_4}{\gamma}\ln\mid a_7
\mid +a_7}{a_2}$ and $\epsilon_6=\frac{a_4}{a_2}$ to make
$\tilde{a_3}=\tilde{a_4}=0$ and $\tilde{a_7}= \pm 1$. Then the
conjugacy class after consider appropriate scaling is $\langle
\alpha X_2 +\beta X_5 + X_7 \rangle$, where $\alpha \neq 0 , \beta \in \mathbb{R}$.\\
\uline{\textbf{Case 6}} $a_1=0$, $a_2 \neq 0$, $a_6= 0$, $a_8= 0$,
$a_7= 0$ : Let $\epsilon_1=\epsilon_2=0$
$\epsilon_5=\frac{a_3}{a_2}$ and $\epsilon_6=\frac{a_4}{a_2}$ to
make $\tilde{a_3} =\tilde{a_4}=0 $. Then the
conjugacy class is $\langle X_2 +\alpha X_5 \rangle$, $\alpha \in \mathbb{R}$.\\
\uline{\textbf{Case 7}} $a_1=0$, $a_2=0$, $a_6 \neq 0$ : Let
$\epsilon_1=-\frac{a_5}{a_6}$, $\epsilon_2=-\frac{a_4}{a_6}$ to
make $\tilde{a_4}=\tilde{a_5}=0$. Then the conjugacy class is of
the form $\langle \alpha X_3 + X_6+ \beta X_7+ \gamma X_8
\rangle$, $\alpha,\beta, \gamma \in \mathbb{R}$.\\
\uline{\textbf{Case 8}} $a_1=0$, $a_2=0$, $a_6 =0$, $a_5 \neq 0$,
$a_8 \neq 0$ : Let $\epsilon_1=-\frac{a_7}{a_8}$ and
$\epsilon_2=\frac{a_4 a_7 - a_3 a_8}{a_5 a_8}$, to get
$\tilde{a_3}=\tilde{a_7}=0$. Then the conjugacy class is of the
form $\langle \alpha X_4+X_5+\beta
X_8 \rangle$, $\alpha \in \mathbb{R}, \beta \neq 0$.\\
\uline{\textbf{Case 9}} $a_1=0$, $a_2=0$, $a_6 =0$, $a_5 \neq 0$,
$a_8=0$, $a_7 \neq 0$ : Let $\epsilon_1=\frac{\ln|a_7|}{\gamma}$
with $\epsilon_2=-\frac{a_4}{\gamma a_5} \ln|a_7|-\frac{a_3}{a_7}$
to make $\tilde{a_3}=0$, $\tilde{a_7}= \pm 1$. Then conjugacy
class after appropriate scaling is of the form $\langle \alpha
X_4+ \beta X_5+ X_7 \rangle$,
$\alpha \in \mathbb{R},\beta \neq 0$.\\
\uline{\textbf{Case 10}} $a_1=0$, $a_2=0$, $a_6 =0$, $a_5 \neq 0$,
$a_8=0$, $a_7= 0$ : Let $\epsilon_2=-\frac{a_3}{a_5}$ to make
$\tilde{a_3}=0$ and so we have the conjugacy class of the form
$\langle \alpha X_4+ X_5\rangle$, with $\alpha \in \mathbb{R}$.\\
\uline{\textbf{Case 11}} $a_1=0$, $a_2=0$, $a_6 =0$, $a_5 =0$,
$a_4 \neq 0$ : Let $\epsilon_1=-\frac{a_3}{a_4}$ to make
$\tilde{a_3}=0$, Then the conjugacy class is $\langle X_4+ \alpha
X_7 + \beta X_8 \rangle$, $\alpha,\beta \in \mathbb{R}$.\\
\uline{\textbf{Case 12}} $a_1=0$, $a_2=0$, $a_6 =0$, $a_5 =0$,
$a_4= 0$, $a_8 \neq 0$ : Let $\epsilon_1=-\frac{a_7}{a_8}$ to have
$\tilde{a_7}=0$. Then the conjugacy class is $\langle \alpha X_3 +
X_8 \rangle$, $\alpha \in \mathbb{R} $.\\
\uline{\textbf{Case 13}} $a_1=0$, $a_2=0$, $a_6 =0$, $a_5 =0$,
$a_4= 0$, $a_8= 0$, $a_7 \neq 0$ : Let
$\epsilon_1=\frac{\ln|a_7|}{\gamma}$ to have $\tilde{a_7}= \pm 1$
and so the conjugacy class after appropriate scaling is
$\langle \alpha X_3+ X_7 \rangle$, $\alpha \in \mathbb{R}$.   \\
\uline{\textbf{Case 14}} $a_1=0$, $a_2=0$, $a_6 =0$, $a_5 =0$,
$a_4= 0$, $a_8= 0$, $a_7= 0$ : Then directly we get a conjugacy
class of the form $\langle X_3 \rangle$.\\
\end{scriptsize}
\end{proof}
\end{theorem}
\subsection {Optimal system for the case $d^2-4 k\rho_2 =\lambda^2$}

The non-zero commutators of the Lie algebra $\mathcal{L}^8$ with
basis (\ref{22}) are given by
\begin{eqnarray} \label{24}
\begin{tabular}{llllll}
$ [X_1,X_4]=X_3 $,&$ [X_1,X_6]=X_5 $ ,& $[X_1,X_7]= -\frac{d}{2 \rho_2}X_7+\frac{\lambda}{2 \rho_2}X_8 $,\\
$[X_1,X_8]=-\frac{d}{2\rho_2}X_8+\frac{\lambda}{2\rho_2}X_7$,&$[X_2,X_5]=X_3 $,& $[X_2,X_6]=X_4$.\\
\end{tabular}
\end{eqnarray}
The Lie algebra $\mathcal{L}^8$ is solvable and the Killing form
is given by $K=2(\hat{a}^2+\hat{b}^2) a^2_1$ where
$\hat{a}=\frac{d}{2 \rho_2}$ and $\hat{b}=\frac{\lambda}{2
\rho_2}$. The adjoint table is given by

\begin{table} [H]  \label{tab:1}
\begin{center}
\begin{tiny}
 \caption {}
\begin{tabular}{ccccccccc}
  \hline
 \rowcolor[gray] {0.9} $Ad (e^\epsilon)$ & $X_1$ & $X_2$ & $X_3$ & $X_4$ & $X_5$ & $X_6$ & $X_7$ & $X_8$ \\
  \hline \hline
  $X_1$ & $X_1$ & $X_2$ & $X_3$ & $X_4-\epsilon X_3$ & $X_5$ & $X_6-\epsilon X_5$ & $ Y_1$ & $Y_2$ \\\hline
  $X_2$ & $X_1$ & $X_2$ & $X_3$ & $X_4$ & $X_5-\epsilon X_3$ & $X_6-\epsilon X_4$ & $X_7$ & $X_8$ \\\hline
  $X_3$ & $X_1$ & $X_2$ & $X_3$ & $X_4$ & $X_5$ & $X_6$ & $X_7$ & $X_8$ \\\hline
  $X_4$ & $X_1+\epsilon X_3$ & $X_2$ & $X_3$ & $X_4$ & $X_5$ & $X_6$ & $X_7$ & $X_8$ \\\hline
  $X_5$ & $X_1$ & $X_2+\epsilon X_3$ & $X_3$ & $X_4$ & $X_5$ & $X_6$ & $X_7$ & $X_8$ \\\hline
  $X_6$ & $X_1+\epsilon X_5$ & $X_2+\epsilon X_4$ & $X_3$ & $X_4$ & $X_5$ & $X_6$ & $X_7$ & $X_8$ \\\hline
  $X_7$ & $X_1-  \epsilon \hat{a} X_7+\epsilon \hat{b}  X_8$ & $X_2$ & $X_3$ & $X_4$ & $X_5$ & $X_6$ & $X_7$ & $X_8$ \\\hline
  $X_8$ & $X_1+ \epsilon \hat{b}  X_7-\epsilon \hat{a}  X_8$ & $X_2$ & $X_3$ & $X_4$ & $X_5$ & $X_6$ & $X_7$ & $X_8$ \\  \hline
\end{tabular}
$Y_1=\frac{1}{2}(e^{\epsilon
(\hat{a}+\hat{b})}+e^{\epsilon(\hat{a}-\hat{b})})X_7+\frac{1}{2}(e^{\epsilon(\hat{a}+\hat{b})}-e^{\epsilon
(\hat{a}-\hat{b})})X_8$ \text{and
}$Y_2=\frac{1}{2}(e^{\epsilon(\hat{a}+\hat{b})}-e^{\epsilon(\hat{a}-\hat{b})})X_7+\frac{1}{2}(e^{\epsilon(\hat{a}+\hat{b})}+e^{\epsilon(\hat{a}-\hat{b})})X_8.$
\end{tiny}
\end{center}
\end{table}
The adjoint group is defined by the matrix
\begin{equation}\label{3-6}
    \begin{array}{c}
      A=Ad(e^{-\epsilon_8 X_8}). Ad(e^{-\epsilon_7 X_7}).
Ad(e^{-\epsilon_6 X_6}). Ad(e^{-\epsilon_5 X_5}).
Ad(e^{-\epsilon_4 X_4}). Ad(e^{-\epsilon_3 X_3}).
Ad(e^{-\epsilon_2 X_2}).
Ad(e^{-\epsilon_1 X_1}), \nonumber \\
    \end{array}
\end{equation}
which is given by
\begin{scriptsize}
\begin{center}
\begin{equation} \label{3-6-1}
A= \begin{pmatrix}
\begin{array}{cccccccc}
  1 & 0 & -\epsilon_4 & 0 & -\epsilon_6 & 0 & \hat{a} \epsilon_7-\hat{b} \epsilon_8 & \hat{a} \epsilon_8-\hat{b} \epsilon_7 \\
  0 & 1 &  -\epsilon_5 & -\epsilon_6 & 0 & 0 & 0 & 0 \\
  0 & 0  & 1 & 0  & 0 & 0 & 0 & 0 \\
  0 & 0 & \epsilon_1 & 1 & 0 & 0 & 0 & 0 \\
  0& 0 & \epsilon_2  & 0 & 1 & 0 & 0 & 0 \\
  0 & 0 & \epsilon_1 \epsilon_2 & \epsilon_2 & \epsilon_1 & 1 & 0 & 0 \\
  0 &0 &  0 & 0 & 0 & 0 & \hat{Y}_1 & \hat{Y}_2 \\
  0 & 0 & 0 & 0 & 0 & 0 &\hat{Y}_2 & \hat{Y}_1 \\
\end{array}
\end{pmatrix}
\end{equation}
$\hat{Y}_1=\frac{1}{2} (e^{- \epsilon_1(\hat{a}-\hat{b})}+e^{-
\epsilon_1(\hat{a}+\hat{b})})$ and $\hat{Y}_2=\frac{1}{2} (e^{-
\epsilon_1(\hat{a}-\hat{b})}-e^{- \epsilon_1(\hat{a}+\hat{b})})$.
\end{center}
\end{scriptsize}
\begin{theorem} An optimal system of one-dimensional Lie algebra
of $L^8$ with basis (\ref{22}) is provided by the following
generators
\begin{scriptsize}
\begin{equation}   \label{29}
\begin{array}{ll|ll}
  X^1=X_1 + \alpha X_2 + \beta X_6, & \alpha \in \mathbb{R},\beta \neq 0, & X^{11}=\alpha X_4 + X_5+ \beta X_8, & \alpha \in \mathbb{R}, \beta \neq 0, \\
  X^2=X_1 + \alpha X_2+ \beta X_4 , & \alpha,\beta \in \mathbb{R}, & X^{12}=\alpha X_4 + \beta X_5 + X_7+X_8, & \alpha \in \mathbb{R},\beta \neq 0,  \\
  X^3=X_2 + \alpha X_6 + \beta X_7+\gamma X_8, & \alpha \neq 0,\beta,\gamma \in \mathbb{R}, &  X^{13}= \alpha X_4+\beta X_5 +X_7 - X_8, & \alpha \in \mathbb{R}, \beta \neq 0, \\
  X^4=X_2 + \alpha X_5 + \beta X_7, & \alpha \in \mathbb{R}, \beta \neq 0, &  X^{14}= \alpha X_4 + X_5, & \alpha \in \mathbb{R}, \\
  X^5=X_2 + \alpha X_5 + \beta X_8, & \alpha \in \mathbb{R}, \beta \neq 0, & X^{15}=  X_4 + \alpha X_7 + \beta X_8, & \alpha,\beta \in \mathbb{R},  \\
  X^6=\alpha X_2 + \beta X_5 + X_7 + X_8, & \alpha \neq 0, \beta \in \mathbb{R}, & X^{16}= \alpha X_3 + X_7, & \alpha \in \mathbb{R}, \\
  X^7=\alpha X_2 + \beta X_5 + X_7 - X_8, & \alpha \neq 0, \beta \in \mathbb{R}, & X^{17}= \alpha X_3 + X_8, & \alpha \in \mathbb{R}, \\
  X^8= X_2 + \alpha X_5, & \alpha \in \mathbb{R}, & X^{18}= \alpha X_3 + X_7 + X_8, & \alpha \in \mathbb{R}, \\
  X^9= \alpha X_3 + X_6+ \beta X_7 + \gamma X_8, & \alpha,\beta,\gamma \in \mathbb{R},& X^{19}= \alpha X_3 + X_7- X_8, & \alpha \in \mathbb{R}, \\
  X^{10}= \alpha X_4 + X_5 + \beta X_7, & \alpha \in \mathbb{R}, \beta \neq 0,& X^{20}=  X_3.&  \\
\end{array}
\end{equation}
\end{scriptsize}

\begin{proof}Let $X$ and $\tilde{X}$ be two elements in the Lie algebra
$\mathcal{L}^8$ with basis (\ref{22}) given by $X=\sum_{i=1}^{8}
a_i X_i$ and $ \tilde{X}=\sum_{i=1}^{8} \tilde{a}_i X_i$. For
simplicity, we will write $X$ and $\tilde{X}$ as row vectors of
the coefficients on the form $X=(\,a_1 \,\, a_2 \,\, . \,\, . \,\,
. \,\, a_8\,)$ and $\tilde{X}=(\,\tilde{a_1} \,\, \tilde{a_2}
\,\, . \,\, . \,\, . \,\, \tilde{a_8}\,)$.\\
Then in order for  $X$ and $\tilde{X}$ to be in the same conjugacy
class, we must have $\tilde{X}=XA$, where $A$ is given by
(\ref{3-6-1}). So, the theorem is proved by solving the system
\begin{equation}\label{3-8-1}
   \begin{array}{l}
   \tilde{a_1}= a_1,\\
   \tilde{a_2}= a_2,\\
   \tilde{a_3}= a_3+\epsilon_1 a_4+\epsilon_2 a_5-\epsilon_5 a_2 -\epsilon_4 a_1+\epsilon_1 \epsilon_2 a_6,\\
   \tilde{a_4}= a_4-\epsilon_6 a_2 +\epsilon_2 a_6, \\
   \tilde{a_5}= a_5-\epsilon_6 a_1+\epsilon_1a_6, \\
   \tilde{a_6}= a_6, \\
   \tilde{a_7}= a_1(\hat{a} \epsilon_7-\hat{b} \epsilon_8)+ \frac{1}{2} a_7  (e^{- \epsilon_1(\hat{a}-\hat{b})}+e^{- \epsilon_1(\hat{a}+\hat{b})})+ \frac{1}{2} a_8  (e^{- \epsilon_1(\hat{a}-\hat{b})}-e^{- \epsilon_1(\hat{a}+\hat{b})}), \\
   \tilde{a_8}= a_1(\hat{a} \epsilon_8-\hat{b} \epsilon_7)+ \frac{1}{2} a_7  (e^{- \epsilon_1(\hat{a}-\hat{b})}-e^{- \epsilon_1(\hat{a}+\hat{b})})+ \frac{1}{2} a_8  (e^{- \epsilon_1(\hat{a}-\hat{b})}+e^{- \epsilon_1(\hat{a}+\hat{b})}), \\
   \end{array}
\end{equation}

for $\{\epsilon_i\}_{i=1}^8$ in term of $\{a_i\}_{i=1}^8$ in order
to get the simplest values of $\{\tilde{a}_i\}_{i=1}^8$. The
results are presented for different cases in the tree diagram (2)
given in the appendix where it is initiated by the sign of the
Killing form and its leafs are given completely.
The full details for each leaf are gives as follows: \\
\begin{scriptsize}
\uline{\textbf{Case 1}} $a_1 \neq 0$, $a_6 \neq 0$ : Let
$\epsilon_1=0$, $\epsilon_2=\frac{a_2 a_5 - a_1 a_4}{a_1 a_6}$,
$\epsilon_4=\frac {a_1 a_3 a_6 - a_1 a_4 a_5+ a_2 a_5^2}{a_1^2
a_6} $, $\epsilon_5=0$, $\epsilon_6=\frac{a_5}{a_1}$,
$\epsilon_7=-\frac{\hat{a}a_7+\hat{b}a_8}{(\hat{a}-\hat{b})a_1}$
and $\epsilon_8 =-\frac{\hat{b}a_7 +\hat{a}a_8} {(\hat{a}-
\hat{b})a_1} $ to have $\tilde{a_3}=\tilde{a_4}=\tilde{a_5}
=\tilde{a_7} =\tilde{a_8}=0$, then we obtain the conjugacy class
$\langle X_1+\alpha X_2 +\beta X_6 \rangle$, where $\alpha \in
\mathbb{R}$, $\beta \neq 0$.\\
\uline{\textbf{Case 2}} $a_1 \neq 0$, $a_6=0$ : Let
$\epsilon_1=\epsilon_2=0$, $\epsilon_4=\frac{a_3}{a_1}$,
$\epsilon_5=0$, $\epsilon_6=\frac{a_5}{a_1}$, $\epsilon_7 =-\frac{
\hat{a}a_7+ \hat{b} a_8}{(\hat{a}-\hat{b})a_1}$ and $\epsilon_8
=-\frac{\hat{b}a_7 +\hat{a}a_8} {(\hat{a}- \hat{b})a_1} $ to have
$\tilde{a_3}=\tilde{a_5}=\tilde{a_6} =\tilde{a_7} =\tilde{a_8}=0$
and so we obtain the donjugacy class $\langle X_1+\alpha X_2
+\beta X_4 \rangle$, where
$\alpha,\beta \in \mathbb{R}$.\\
\uline{\textbf{Case 3}} $a_1 = 0$, $a_2\neq 0$, $a_6\neq0$ : Let
$\epsilon_1=-\frac{a_5}{a_6}$, $\epsilon_2=0$, $\epsilon_5= \frac
{a_3 a_6 - a_4 a_5}{a_2 a_6}$ and $\epsilon_6=\frac{a_4}{a_2}$ to
have $\tilde{a_3}=\tilde{a_4}=\tilde{a_5}=0$. Then the conjugacy
class is of the form $\langle X_2+ \alpha X_6 +\beta X_7 + \gamma
X_8 \rangle$, where $\alpha \neq 0$, $\beta ,
\gamma \in \mathbb{R}$.\\
\uline{\textbf{Case 4}} $a_1 = 0$, $a_2 \neq 0$, $a_6=0$, $ a_7 +
a_8 \neq 0 $, $\frac{a_7-a_8}{a_7+a_8}>0$  : Let
$\epsilon_1=\frac{1}{2 \hat{b}} \ln (\frac{a_7-a_8}{a_7+a_8})$,
$\epsilon_2=0$ and $\epsilon_5=\frac{a_4 }{2 \hat{b}a_2} \ln
(\frac{a_7-a_8}{a_7+a_8}) +\frac{a_3}{a_2}$ and
$\epsilon_6=\frac{a_4}{a_2}$ to have
$\tilde{a_3}=\tilde{a_4}=\tilde{a_8}=0$ and so $\langle X_2
+\alpha X_5 + \beta X_7 \rangle$,
where $\alpha \in \mathbb{R}$, $\beta \neq 0$.\\
\uline{\textbf{Case 5}} $a_1 = 0$, $a_2 \neq 0$, $a_6=0$, $ a_7 +
a_8 \neq 0 $, $\frac{a_7-a_8}{a_7+a_8}< 0$  : Let
$\epsilon_1=\frac{1}{2 \hat{b}} \ln (-\frac{a_7-a_8}{a_7+a_8})$,
$\epsilon_2=0$ and $\epsilon_5=\frac{a_4 }{2 \hat{b}a_2} \ln
(-\frac{a_7-a_8}{a_7+a_8}) +\frac{a_3}{a_2}$, and
$\epsilon_6=\frac{a_4}{a_2}$ to have
$\tilde{a_3}=\tilde{a_4}=\tilde{a_7}=0$ and so we obtain the
conjugacy class $\langle X_2
+\alpha X_5 + \beta X_8 \rangle$, where $\alpha \in \mathbb{R} $, $\beta \neq 0$.\\
\uline{\textbf{Case 6}} $a_1 = 0$, $a_2 \neq 0$, $a_6=0$, $ a_7 +
a_8 \neq 0 $, $\frac{a_7-a_8}{a_7+a_8}= 0$ : Note that in this
case $a_8 \neq 0$, so let $\epsilon_1=\frac{\ln|a_8|}{\hat{a}-
\hat{b}}$, $\epsilon_2=0$,
$\epsilon_5=\frac{a_4}{(\hat{a}-\hat{b})a_2} \ln|a_8|
+\frac{a_3}{a_2}$ and $\epsilon_6=\frac{a_4}{a_2}$ to have
$\tilde{a_3}=\tilde{a_4}=0$ and so $\langle \alpha X_2 + \beta X_5
+ X_7 + X_8 \rangle$, where $\alpha \neq 0 $
, $\beta \in \mathbb{R} $.\\
\uline{\textbf{Case 7}} $a_1 = 0$, $a_2 \neq 0$, $a_6=0$, $ a_7 +
a_8= 0 $, $ a_8 \neq 0$ : Let $\epsilon_1=\frac{\ln|a_8|}{\hat{a}+
\hat{b}}$, $\epsilon_2=0$, $\epsilon_5= \frac{a_4}{(\hat{a}
+\hat{b})a_2} \ln|a_8| +\frac{a_3}{a_2}$ and $\epsilon_6=\frac{a_4
}{a_2}$ to have $\tilde{a_3}=\tilde{a_4}=0$ and so we obtain the
congugacy class $\langle \alpha X_2 + \beta X_5 + X_7 - X_8
\rangle$, where $\alpha \neq 0 $
, $\beta \in \mathbb{R} $.\\
\uline{\textbf{Case 8}} $a_1 = 0$, $a_2 \neq 0$, $a_6=0$, $ a_7 +
a_8= 0 $, $ a_8 = 0$ : Let $\epsilon_1=\epsilon_2=0$, $\epsilon_5
=\frac{a_3}{a_2}$ and $\epsilon_6= \frac{a_4}{a_2}$ to have
$\tilde{a_3}=\tilde{a_4}=0$ and so the conjugacy class is of the
form $\langle X_2 + \alpha X_5
\rangle $ , where $\alpha \in \mathbb{R} $.\\
\uline{\textbf{Case 9}} $a_1 = 0$, $a_2= 0$, $a_6 \neq0$ : Let
$\epsilon_1=-\frac{a_5}{a_6} $,  $\epsilon_2=-\frac{a_4}{a_6}$ and
$\epsilon_6=0$ to have $\tilde{a_4}=\tilde{a_5}=0$ and so the
conjugacy class is of the form $\langle \alpha X_3 + X_6 +\beta
X_7 + \gamma X_8
\rangle $ , where $\alpha,\beta,\gamma \in \mathbb{R} $.\\
\uline{\textbf{Case 10}} $a_1 = 0$, $a_2=0$, $a_6=0$, $a_5 \neq
0$, $ a_7 + a_8 \neq 0 $, $\frac{a_7-a_8}{a_7+a_8}>0$  : Let
$\epsilon_1=\frac{1}{2 \hat{b}} \ln (\frac{a_7-a_8}{a_7+a_8})$ and
$\epsilon_2=\frac{a_4 }{2 \hat{b}a_5} \ln
(\frac{a_7-a_8}{a_7+a_8}) +\frac{a_3}{a_5}$ to have
$\tilde{a_3}=\tilde{a_8}=0$ and so the conjugacy class is of the
form $\langle \alpha X_4 + X_5 +
\beta X_7 \rangle$, where $\alpha \in \mathbb{R}$, $\beta \neq 0$.\\
\uline{\textbf{Case 11}} $a_1 = 0$, $a_2=0$, $a_6=0$, $a_5 \neq
0$, $ a_7 + a_8 \neq 0 $, $\frac{a_7-a_8}{a_7+a_8}<0$  : Let
$\epsilon_1=\frac{1}{2 \hat{b}} \ln (-\frac{a_7-a_8}{a_7+a_8})$
and $\epsilon_2=-\frac{a_4 }{2 \hat{b}a_5} \ln
(-\frac{a_7-a_8}{a_7+a_8}) - \frac{a_3}{a_5}$ to have
$\tilde{a_3}=\tilde{a_7}=0$ and so the conjugacy class is of the
form $\langle \alpha X_4 + X_5 +
\beta X_8 \rangle$, where $\alpha \in \mathbb{R}$, $\beta \neq 0$.\\
\uline{\textbf{Case 12}} $a_1 = 0$, $a_2=0$, $a_6=0$, $a_5 \neq
0$, $ a_7 + a_8 \neq 0 $, $\frac{a_7-a_8}{a_7+a_8}=0$ : Note that
in this case $a_8 \neq 0$, so let
$\epsilon_1=\frac{\ln|a_8|}{\hat{a} - \hat{b}}$ and
$\epsilon_2=-\frac{a_4 }{(\hat{a} - \hat{b})a_5} \ln|a_8| -
\frac{a_3}{a_5}$, to have $\tilde{a_3}=0$, and so the conjugacy
class is $\langle \alpha X_4 + \beta X_5 +
X_7+ X_8 \rangle$, where $\alpha \in \mathbb{R}$, $\beta \neq 0$.\\
\uline{\textbf{Case 13}} $a_1 = 0$, $a_2=0$, $a_6=0$, $a_5 \neq
0$, $ a_7 + a_8= 0 $, $ a_8 \neq 0$ : Let
$\epsilon_1=\frac{\ln|a_8|}{\hat{a}+ \hat{b}}$ and $\epsilon_2=
\frac{a_4}{(\hat{a} +\hat{b})a_5} \ln|a_8| +\frac{a_3}{a_5}$ to
have $\tilde{a_3}=0$ and so the conjugacy class is of the form
$\langle \alpha X_4 + \beta X_5 + X_7
- X_8 \rangle$, where $\alpha \in \mathbb{R}$, $\beta \neq 0 $.\\
\uline{\textbf{Case 14}} $a_1 = 0$, $a_2=0$, $a_6=0$, $a_5 \neq
0$, $ a_7 + a_8= 0 $, $ a_8 = 0$ : Let $\epsilon_1=0$ and
$\epsilon_2 =-\frac{a_3}{a_5}$ to have $\tilde{a_3}=0$ and so we
obtain the conjugacy class $\langle \alpha
X_4 + X_5 \rangle $, $\alpha \in \mathbb{R} $.\\
\uline{\textbf{Case 15}} $a_1=0$, $a_2=0$, $a_6=0$, $a_5=0$, $a_4
\neq 0$ : let $\epsilon_1=-\frac{a_3}{a_4}$ to have
$\tilde{a_3}=0$ and so the conjugacy class is $\langle X_4+\alpha
X_7 +\beta X_8 \rangle$, where $\alpha,\beta \in
\mathbb{R}$.\\
\uline{\textbf{Case 16}} $a_1 = 0$, $a_2=0$, $a_6=0$, $a_5=0$,
$a_4=0$, $ a_7 + a_8 \neq 0 $, $\frac{a_7-a_8}{a_7+a_8}>0$ : Let
$\epsilon_1=\frac{1}{2 \hat{b}} \ln (\frac{a_7-a_8}{a_7+a_8})$ to
have $\tilde{a_8}=0$ and so the conjugacy class after appropriate
scaling is of the form
$\langle \alpha X_3 + X_7\rangle$, where $\alpha \in \mathbb{R}$.\\
\uline{\textbf{Case 17}} $a_1 = 0$, $a_2=0$, $a_6=0$, $a_5=0$,
$a_4=0$, $ a_7 + a_8 \neq 0 $, $\frac{a_7-a_8}{a_7+a_8}<0$ : Let
$\epsilon_1=\frac{1}{2 \hat{b}} \ln (-\frac{a_7-a_8}{a_7+a_8})$ to
have $\tilde{a_7}=0$ and so the conjugacy class after appropriate
scaling is of the form
$\langle \alpha X_3 + X_8 \rangle$, where $\alpha \in \mathbb{R}$.\\
\uline{\textbf{Case 18}} $a_1 = 0$, $a_2=0$, $a_6=0$, $a_5=0$,
$a_4=0$, $ a_7 + a_8 \neq 0 $, $\frac{a_7-a_8}{a_7+a_8}=0$ : Note
that $a_8 \neq 0$ so let $\epsilon_1=\frac{\ln|a_8|}{\hat{a} -
\hat{b}}$ to have thr conjugacy class
$\langle \alpha X_3 +X_7 + X_8 \rangle$, where $\alpha \in \mathbb{R}$.\\
\uline{\textbf{Case 19}} $a_1 = 0$, $a_2=0$, $a_6=0$, $a_5=0$,
$a_4=0$, $ a_7 + a_8= 0 $, $a_8 \neq 0 $ : Let
$\epsilon_1=\frac{\ln|a_8|}{\hat{a} + \hat{b}}$ to have the
conjugacy class $\langle \alpha X_3 +X_7 - X_8 \rangle$, where $\alpha \in \mathbb{R}$.\\
\uline{\textbf{Case 20}} $a_1 = 0$, $a_2=0$, $a_6=0$, $a_5=0$,
$a_4=0$, $ a_7 + a_8= 0 $, $a_8 = 0 $ : Then directly we have the
conjugacy
class $\langle X_3  \rangle$.\\
\end{scriptsize}
\end{proof}
\end{theorem}

\subsection {Optimal system for the case $d^2-4 k\rho_2 = -\mu^2$}

The non-zero commutators of the Lie algebra $\mathcal{L}^8$ with
basis (\ref{25}) are given by
\begin{eqnarray} \label{27}
\begin{tabular}{llllll}
$ [X_1,X_4]=X_3 $,& $ [X_1,X_6]=X_5 $ ,&$[X_1,X_7]= -\frac{d}{2 \rho_2}X_7-\frac{\mu}{2 \rho_2}X_8 $, \\
$[X_1,X_8]=-\frac{d}{2\rho_2}X_8+\frac{\mu}{2\rho_2}X_7$,&$[X_2,X_5]=X_3 $,& $[X_2,X_6]=X_4$.\\
\end{tabular}
\end{eqnarray}
The Lie algebra $\mathcal{L}^8$ is solvable and the Killing form
is given by $K=2(\hat{a}^2-\hat{b}^2) a^2_1$ where
$\hat{a}=\frac{d}{2\rho_2}$ and $\hat{b}=\frac{\mu}{2\rho_2}$. The
adjoint table is given by

\begin{table} [H]  \label{tab:1}
\begin{center}
\begin{tiny}
 \caption {}
\begin{tabular}{ccccccccc}
  \hline
  \rowcolor[gray] {0.9} $Ad (e^\epsilon)$ & $X_1$ & $X_2$ & $X_3$ & $X_4$ & $X_5$ & $X_6$ & $X_7$ & $X_8$ \\ \hline  \hline
  $X_1$ & $X_1$ & $X_2$ & $X_3$ & $X_4-\epsilon X_3$ & $X_5$ & $X_6-\epsilon X_5$ & $ Y_1$ & $Y_2$ \\\hline
  $X_2$ & $X_1$ & $X_2$ & $X_3$ & $X_4$ & $X_5-\epsilon X_3$ & $X_6-\epsilon X_4$ & $X_7$ & $X_8$ \\\hline
  $X_3$ & $X_2$ & $X_2$ & $X_3$ & $X_4$ & $X_5$ & $X_6$ & $X_7$ & $X_8$ \\\hline
  $X_4$ & $X_1+\epsilon X_3$ & $X_2$ & $X_3$ & $X_4$ & $X_5$ & $X_6$ & $X_7$ & $X_8$ \\\hline
  $X_5$ & $X_1$ & $X_2+\epsilon X_3$ & $X_3$ & $X_4$ & $X_5$ & $X_6$ & $X_7$ & $X_8$ \\\hline
  $X_6$ & $X_1+\epsilon X_5$ & $X_2+\epsilon X_4$ & $X_3$ & $X_4$ & $X_5$ & $X_6$ & $X_7$ & $X_8$ \\\hline
  $X_7$ & $X_1-  \epsilon \hat{a} X_7-\epsilon \hat{b}  X_8$ & $X_2$ & $X_3$ & $X_4$ & $X_5$ & $X_6$ & $X_7$ & $X_8$ \\\hline
  $X_8$ & $X_1+ \epsilon \hat{b}  X_7-\epsilon \hat{a}  X_8$ & $X_2$ & $X_3$ & $X_4$ & $X_5$ & $X_6$ & $X_7$ & $X_8$ \\  \hline
\end{tabular}\\
$Y_1=e^{\epsilon \hat{a}}(\cos(\epsilon \hat{b})X_7+\sin(\epsilon
\hat{b})X_8) $ \text{and } $Y_2=e^{\epsilon
\hat{a}}(-\sin(\epsilon \hat{b})X_7+\cos(\epsilon \hat{b})X_8).$
\end{tiny}
\end{center}
\end{table}
The adjoint group is defined by the matrix
\begin{equation}\label{}
    \begin{array}{c}
      A=Ad(e^{-\epsilon_8 X_8}). Ad(e^{-\epsilon_7 X_7}).
Ad(e^{-\epsilon_6 X_6}). Ad(e^{-\epsilon_5 X_5}).
Ad(e^{-\epsilon_4 X_4}). Ad(e^{-\epsilon_3 X_3}).
Ad(e^{-\epsilon_2 X_2}).
Ad(e^{-\epsilon_1 X_1}), \nonumber \\
    \end{array}
\end{equation}
which is given by
\begin{scriptsize}
\begin{equation} \label{3-9-1}
 A= \begin{pmatrix}
\begin{array}{cccccccc}
  1 & 0 & -\epsilon_4 & 0 & -\epsilon_6 & 0 & \hat{a} \epsilon_7-\hat{b} \epsilon_8 & \hat{a} \epsilon_8+\hat{b} \epsilon_7 \\
  0 & 1 & -\epsilon_5 & -\epsilon_6 & 0 & 0 & 0 & 0 \\
  0 & 0 & 1 & 0 & 0 & 0 & 0 & 0 \\
  0 & 0 & \epsilon_1 & 1 & 0 & 0 & 0 & 0 \\
  0 & 0 & \epsilon_2 & 0 & 1 & 0 & 0 & 0 \\
  0 & 0 &\epsilon_1 \epsilon_2 & \epsilon_2 & \epsilon_1 & 1 & 0 & 0 \\
  0 & 0 & 0 & 0 & 0 & 0 & \hat{Y}_1 & -\hat{Y}_2 \\
  0 & 0 & 0 & 0 & 0 & 0 & \hat{Y}_2 & \hat{Y}_1 \\
\end{array}
 \end{pmatrix}.
\end{equation}
\end{scriptsize}
 \begin{scriptsize} \begin{equation}\label{3-10}
    \begin{array}{cc}
      \hat{Y}_1=e^{-\epsilon_1 \hat{a}} \cos(\epsilon_1 \hat{b}), & \hat{Y}_2=e^{-\epsilon_1 \hat{a}} \sin(\epsilon_1 \hat{b}) \nonumber \\
    \end{array}
\end{equation}\end{scriptsize}

\begin{theorem} An optimal system of one-dimensional Lie algebra
$L^8$ with basis (\ref{25}) is provided by the following
generators
\begin{scriptsize}
\begin{equation} \label {30}
\begin{array}{ll|ll}
  X^1=X_1 + \alpha X_2 + \beta X_6, & \alpha \in \mathbb{R},\beta \neq 0, & X^6=X_4 + \alpha X_5 + \beta X_7+\gamma X_8, & \alpha, \beta , \gamma \in \mathbb{R}, \\
  X^2=X_1 + \alpha X_2 + \beta X_4 , & \alpha, \beta \in \mathbb{R},  & X^7=X_5 + \alpha X_8, & \alpha \neq 0, \\
  X^3=X_2 + \alpha X_6 + \beta X_7+\gamma X_8, & \alpha \neq 0,\beta , \gamma \in \mathbb{R}, & X^8= X_5 + \alpha X_7, & \alpha \in \mathbb{R}, \\
  X^4=X_2 + \alpha X_5 + \beta X_7+\gamma X_8, & \alpha , \beta , \gamma \in \mathbb{R}, &  X^9=\alpha X_3 + X_8, & \alpha \in \mathbb{R},\\
  X^5=\alpha X_3 +  X_6 + \beta X_7 + \gamma X_8, & \alpha, \beta , \gamma \in \mathbb{R}, & X^{10}=\alpha X_3 + \beta X_7, & \alpha,\beta \in \mathbb{R}. \\
\end{array}
\end{equation}
\end{scriptsize}
\begin{proof}
Let $X$ and $\tilde{X}$ be two elements in the Lie algebra
$\mathcal{L}^8$ with basis (\ref{25}) given by $X=\sum_{i=1}^{8}
a_i X_i$ and $ \tilde{X}=\sum_{i=1}^{8} \tilde{a}_i X_i$. For
simplicity, we will write $X$ and $\tilde{X}$ as row vectors of
the coefficients on the form $X=(\,a_1 \,\, a_2 \,\, . \,\, . \,\,
. \,\, a_8\,)$ and $\tilde{X}=(\,\tilde{a_1} \,\, \tilde{a_2} \,\,
. \,\, . \,\, . \,\, \tilde{a_8}\,)$. Then in order for $X$ and
$\tilde{X}$ to be in the same conjugacy class, we must have
$\tilde{X}=XA$, where $A$ is given by (\ref{3-9-1}). So, the
theorem is proved by solving the system
\begin{equation}\label{3-11}
   \begin{array}{l}
    \tilde{a_1}= a_1, \\
   \tilde{a_2}= a_2, \\
   \tilde{a_3}= a_3+\epsilon_1 a_4+\epsilon_2 a_5-\epsilon_5 a_2 -\epsilon_4 a_1+\epsilon_1 \epsilon_2 a_6,\\
   \tilde{a_4}= a_4-\epsilon_6 a_2 +\epsilon_2 a_6, \\
   \tilde{a_5}= a_5-\epsilon_6 a_1+\epsilon_1a_6, \\
   \tilde{a_6}= a_6, \\
   \tilde{a_7}= a_1(\hat{a} \epsilon_7-\hat{b}\epsilon_8)+ e^{-\epsilon_1 \hat{a}}\left(\cos(\epsilon_1 \hat{b} )a_7+\sin(\epsilon_1 \hat{b} )a_8 \right), \\
   \tilde{a_8}= a_1(\hat{b} \epsilon_7+\hat{a}\epsilon_8)- e^{-\epsilon_1 \hat{a}}\left(\sin(\epsilon_1 \hat{b} )a_7-\cos(\epsilon_1 \hat{b} )a_8 \right), \\
   \end{array}
\end{equation}

for $\{\epsilon_i\}_{i=1}^8$ in term of $\{a_i\}_{i=1}^8$ in order
to get the simplest values of $\{\tilde{a}_i\}_{i=1}^8$. The
results are presented for different cases in the tree diagram (3)
given in the appendix where it is initiated by the sign of the
Killing form and its leafs are given completely.
The full details for each leaf are given as follows: \\
\begin{scriptsize}
\uline{\textbf{Case 1}} $a_1 \neq 0$, $a_6 \neq 0$ : Let
$\epsilon_1=-\frac{a_5}{a_6}$, $\epsilon_2=-\frac{a_4}{a_6}$,
$\epsilon_4=\frac{a_3 a_6 -a_4 a_5}{a_1 a_6}$,
$\epsilon_5=\epsilon_6=0$,\\
$\epsilon_7=\frac{1}{a_1(\hat{a}^2+\hat{b}^2)}e^{\frac{a_5\hat{a}}{a_6}}\left(\left(\hat{a}a_8-\hat{b}a_7\right)\sin\left(\frac{a_5\hat{b}}{a_6}\right)-\left(\hat{a}a_7+\hat{b}a_8\right)\cos\left(\frac{a_5\hat{b}}{a_6}\right)\right)$
and \\
$\epsilon_8=-\frac{1}{a_1(\hat{a}^2+\hat{b}^2)}e^{\frac{a_5\hat{a}}{a_6}}\left(\left(\hat{a}a_7+\hat{b}a_8\right)\sin\left(\frac{a_5\hat{b}}{a_6}\right)+\left(\hat{a}a_8-\hat{b}a_7\right)\cos\left(\frac{a_5\hat{b}}{a_6}\right)\right)$
to have
$\tilde{a_3}=\tilde{a_4}=\tilde{a_5}=\tilde{a_7}=\tilde{a_8}=0$
then we have the conjugacy class $\langle X_1+\alpha X_2 +\beta
X_6 \rangle$, with $\alpha \in \mathbb{R}$, and $\beta \neq
0$.\\
\uline{\textbf{Case 2}} $a_1 \neq 0$, $a_6=0$ : Let
$\epsilon_1=\epsilon_2=0$, $\epsilon_4=\frac{a_3}{a_1}$,
$\epsilon_6=\frac{a_5}{a_1}$, $\epsilon_7=-\frac{\hat{a}
a_7+\hat{b} a_8}{a_1(\hat{a}^2+\hat{b}^2)}$, and
$\epsilon_8=-\frac{\hat{a} a_8-\hat{b}
a_7}{a_1(\hat{a}^2+\hat{b}^2)}$ to have
$\tilde{a_3}=\tilde{a_5}=\tilde{a_7}=\tilde{a_8}=0$. Then we have
the conjugacy class of the form $\langle X_1+\alpha X_2+\beta X_4
\rangle$ where $\alpha,\beta \in \mathbb{R}$.\\
\uline{\textbf{Case 3}} $a_1=0$, $a_2 \neq 0$, $a_6 \neq 0$ : Let
$\epsilon_1=-\frac{a_5}{a_6}$, $\epsilon_2=0$,
$\epsilon_5=\frac{a_3 a_6-a_4 a_5}{a_2 a_6}$, and
$\epsilon_6=\frac{a_4}{a_2}$ to get
$\tilde{a_3}=\tilde{a_4}=\tilde{a_5}=0$.Then the conjugacy class
is of the form $\langle X_2+\alpha X_6+\beta X_7+\gamma X_8
\rangle$, where
$\alpha \neq 0$, $\beta,\gamma \in \mathbb{R}$.\\
\uline{\textbf{Case 4}} $a_1=0$, $a_2 \neq 0$, $a_6=0$ : Let
$\epsilon_1=\epsilon_2=0$, $\epsilon_5=\frac{a_3}{a_2}$,
$\epsilon_6=\frac{a_4}{a_2}$ to get $\tilde{a_3}=\tilde{a_4}=0$.
Then the conjugacy class is $\langle X_2+\alpha X_5+\beta X_7
+\gamma X_8 \rangle$, $\alpha,\beta,\gamma \in
\mathbb{R}$.\\
\uline{\textbf{Case 5}} $a_1=0$, $a_2= 0$, $a_6 \neq 0$ : Let
$\epsilon_1=-\frac{a_5}{a_6}$, $\epsilon_2=-\frac{a_4}{a_6}$ to
have $\tilde{a_4}=\tilde{a_5}=0$ and so we obtain the conjugacy
class $\langle \alpha X_3+X_6+\beta X_7
+\gamma X_8 \rangle$, $\alpha,\beta,\gamma \in \mathbb{R}$.\\
\uline{\textbf{Case 6}} $a_1=0$, $a_2= 0$, $a_6 = 0$, $a_4 \neq 0$
: Let $\epsilon_1=-\frac{a_3}{a_4}$ to have $\tilde{a_3}=0$ and so
we obtain the conjugacy class $\langle X_4 + \alpha
X_5 + \beta X_7 + \gamma X_8 \rangle$, $\alpha,\beta,\gamma \in \mathbb{R}$.\\
\uline{\textbf{Case 7}} $a_1=0$, $a_2= 0$, $a_6 = 0$, $a_4= 0$ $
a_5 \neq 0$, $a_8 \neq 0$ : Let
$\epsilon_1=-\frac{1}{b}\tan^{-1}(\frac{a_7}{a_8})$ and
$\epsilon_2=-\frac{a_3}{a_5}$ to have $\tilde{a_3}=\tilde{a_7}=0$.
Then the conjugacy class is of the form $\langle X_5+\alpha X_8 \rangle$, where $\alpha \neq 0$.\\
\uline{\textbf{Case 8 }} $a_1=0$, $a_2= 0$, $a_6 = 0$, $a_4= 0$ $
a_5 \neq 0$, $a_8 = 0$ : Let $\epsilon_1=0$ and
$\epsilon_2=-\frac{a_3}{a_5}$ to have $\tilde{a_3}=\tilde{a_8}=0$.
Then the conjugacy class is $\langle X_5+\alpha X_7 \rangle$, $\alpha \in \mathbb{R}$.\\
\uline{\textbf{Case 9 }} $a_1=0$, $a_2= 0$, $a_6 = 0$, $a_4= 0$,
$a_5=0$, $a_8 \neq 0$ : Let
$\epsilon_1=-\frac{1}{b}\tan^{-1}(\frac{a_7}{a_8})$ to have
$\tilde{a_7}=0$. By appropriate scaling the conjugacy class is of
the form
 $\langle \alpha X_3+ X_8 \rangle$, where $\alpha \in \mathbb{R}$.\\
\uline{\textbf{Case 10 }} $a_1=0$, $a_2= 0$, $a_6 = 0$, $a_4= 0$,
$a_5=0$, $a_8=0$ : Let $\epsilon_1=0$
to have $\tilde{a_8}=0$ and so we have the conjugacy class $\langle \alpha X_3+\beta X_7 \rangle$, where $\alpha,\beta \in \mathbb{R}$.\\
\end{scriptsize}
\end{proof}
\end{theorem}
\section{Optimal reductions and invariant solutions}
It is known that the invariant solutions for PDEs can be
determined by two procedures, which are the invariant form method
and the direct substitution method \cite{Bluman}. The idea of
looking for group invariant solutions generalize quite naturally
to PDEs with any number of independent and dependent variables. A
one parameter group that acts nontrivially on one or more
independent variables can be used to reduce the number of
independent variables by one.\\
In this section, we focus on the invariant form method which
requires that at least one of the infinitesimals $\xi^1$ and
$\xi^2$ does not equal zero \cite{Bluman,hydon}. Hence, we solve
the invariance surface conditions explicitly by solving the
corresponding characteristic equation given by
\begin{equation}  \label{31}
\frac{d t }{\xi^1 (t,x,\varphi, \psi)} = \frac{d x }{\xi^2
(t,x,\varphi, \psi)}=\frac{d \varphi }{\eta^1 (t,x,\varphi,
\psi)}= \frac{d \psi }{\eta^2 (t,x,\varphi, \psi)},
\end{equation}
to get the corresponding invariants which are used to reduce the
number of independent variables by one. The procedure is explained
in details in three examples. Moreover,  all possible invariant
variables  of the optimal systems (\ref{28}), (\ref{29}) and
(\ref{30}) and their corresponding reductions for the three
non-linear cases of $\chi(\psi_x)$ are given in Table 4, Table 5
and Table 6, respectively.

\begin{example} Reduction for case (2.2.1) using invariant form of $X^3$.\\
Consider the generator $X^3=X_2 +\alpha X_6+\beta X_7 +\gamma X_8$
where $\alpha \neq 0,\beta,\gamma \in \mathbb{R}$, from the
optimal system (\ref{28}). Solving the corresponding
characteristic equation
\begin{equation}  \label{31}
\frac{d t }{0} = \frac{d x }{1}=\frac{d \varphi }{\alpha t x}=
\frac{d \psi }{\alpha(\frac{d}{k}-t)+e^{-\frac{d}{2
\rho_2}t}(\beta+\gamma t)}, \nonumber
\end{equation}
gives the invariant variables as follows:
\begin{equation}\label{e1}
    \begin{array}{ccc}
      \phi (t,x)\, = \,Z(\zeta ) + \frac{\alpha }{2}tx^2, & \psi (t,x)\,
= \,W(\zeta ) + \alpha x(\frac{d}{k} - t) + x(\beta+ \gamma t)e^{
- \frac{d}{{2 \rho_2 }} t}, & \zeta \, = \,t. \\
    \end{array}
\end{equation}
 The reduced system
resulting from the invariant variables (\ref{e1}) is the system of
ODEs of the form
\begin{equation} \label {e2}
\begin{array}{l}
   \rho _1 \,Z'' -k(\gamma \zeta  + \beta )e^{ -  \frac{d}{{2 \rho_2 }} \zeta} {\rm{ }} - \alpha d\, = \,0, \\
    \rho _2 \,W'' + d\,W' + k\,W\, = \,0. \\ \nonumber
\end{array}
\end{equation}
Solving this system yields the following solution
\begin{equation} \label {34}
\begin{array}{l}
  \varphi (t,x) = \frac{{\rho _2 }}{{\rho _1 d}}\left( {d\gamma t +d\beta  + 4\gamma \rho _2 } \right)e^{ - \frac{d}{{2\rho _2 }}t}+ \alpha t x^2    + \frac{{\alpha d}}{{\rho _1 }}t^2 + 2c_1 t +2c_2 ,\\
  \psi (t,x) = (\gamma tx + \beta x + c_3 t + c_4 )e^{ -\frac{{d}}{{2\rho _2 }}t}  - \alpha t x + \frac{\alpha
  d}{k}x,\\ \nonumber
\end{array}
\end{equation}
where $d=2\sqrt{k\rho_2}$.
\end{example}

\begin{example} Reduction for case (2.2.2) using invariant form of $X^3$.\\
Consider the generator $X^3=X_2+\alpha X_6 +\beta X_7+\gamma X_8 $
where $\alpha \neq 0,\beta,\gamma \in \mathbb{R}$, from the
optimal system (\ref{29}). Solving the corresponding
characteristic equation
\begin{equation}  \label{31}
\frac{d t }{0} = \frac{d x }{1}=\frac{d \varphi }{\alpha t x}=
\frac{d \psi }{\alpha(\frac{d}{k}-t)+e^{-\frac{d}{2
\rho_2}t}\left(\beta \cosh\left(\frac{\lambda ~t}{2
\rho_2}\right)+\gamma \sinh\left(\frac{\lambda ~t}{2
\rho_2}\right)\right)}, \nonumber
\end{equation}
gives the invariant variables as follows:
\begin{equation}\label{e11}
    \begin{array}{ccc}
      \phi (t,x)\, = \,Z(\zeta ) + \frac{\alpha }{2}tx^2, & \psi (t,x)\,= \,W(\zeta ) + \alpha x(\frac{d}{k} - t) + xe^{ - \frac{d}{{2\rho_2 }}t} (\beta \cosh (\frac{{\lambda~t}}{{2\rho _2 }}) + \gamma\sinh (\frac{{\lambda~t}}{{2\rho _2 }})), & \zeta \, = \,t.
      \\
    \end{array}
\end{equation}
 The reduced system
resulting from the invariant variables (\ref{e11}) is the system
of ODEs of the form
\begin{equation} \label {e2}
\begin{array}{l}
   \rho _1 \,Z'' - k\,e^{ -\frac{d}{{2\rho _2 }}\zeta }(\beta \cosh (\frac{{\lambda \zeta }}{{2\rho _2}}) + \gamma \sinh (\frac{{\lambda \zeta }}{{2\rho _2 }}))  - \alpha d\, = \,0, \\
    \rho _2 \,W'' + d\,W' + k\,W\, = \,0. \\ \nonumber
\end{array}
\end{equation}
Solving this system yields the following solution
\begin{equation} \label {34}
\begin{array}{ll}
  \varphi (t,x) =&\frac{1}{{2k\rho _1 }}\left( {\left( {2\,\rho _2 \beta \,k + \beta\,\lambda ^2  + \lambda d\gamma } \right)\cosh \left({\,\frac{{\lambda \,t}}{{2\rho _2 }}} \right) + \left( {2\,\gamma\,\rho _2 k + \,\lambda \beta \,d + \gamma \,\lambda ^2 }\right)\sinh \left( {\,\frac{{\lambda \,t}}{{2\rho _2 }}} \right)}\right){\rm{e}}^{ - \,\frac{d}{{2\rho _2 }}t}\\
   &+\frac{{\alpha\,d}}{{2\rho _1 }}{\rm{ }}t^2  + \frac{\alpha }{2}\,t\,x^2  +\,c_1 t + c_2 ,\\
  \psi (t,x) =&  \left( {(c_3+\beta x) \cosh \left({\frac{{\lambda \,t}}{{2\rho _2 }}} \right) + (c_4+\gamma x) \sinh \left({\frac{{\lambda \,t}}{{2\rho _2 }}} \right)} \right){\rm{e}}^{ -\,\frac{d}{{2\rho _2 }}t} - x\alpha \,t+ \frac{{\alpha \,d}}{k}x,\\
  \nonumber
\end{array}
\end{equation}
where $d=\sqrt{k\rho_2+\lambda^2}$.
\end{example}

\begin{example} Reduction for case (2.2.3) using invariant form of $X^3$.\\
Consider the generator $X^3=X_2+\alpha X_6 +\beta X_7+\gamma X_8 $
where $\alpha \neq 0,\beta,\gamma \in \mathbb{R}$, from the
optimal system (\ref{30}). Solving the corresponding
characteristic equation
\begin{equation}  \label{}
\frac{d t }{0} = \frac{d x }{1}=\frac{d \varphi }{\alpha t x}=
\frac{d \psi }{\alpha(\frac{d}{k}-t)+e^{-\frac{d}{2
\rho_2}t}\left(\beta \cos\left(\frac{\mu t}{2
\rho_2}\right)+\gamma \sin\left(\frac{\mu t}{2
\rho_2}\right)\right)}, \nonumber
\end{equation}
gives the invariant variables as follows:
\begin{equation}\label{e111}
    \begin{array}{ccc}
      \phi (t,x)\, = \,Z(\zeta ) + \frac{\alpha }{2}tx^2 ,& \psi(t,x) \, = \,W(\zeta ) + \alpha x(\frac{d}{k} - t) + xe^{- \frac{d}{{2\rho _2 }}t}\left(\beta \cos (\frac{{\mu t}}{{2\rho _2 }}) + \gamma \sin (\frac{{\mu t}}{{2\rho _2 }})\right), & \zeta \, = \,t.
      \\
    \end{array}
\end{equation}
 The reduced system
resulting from the invariant variables (\ref{e111}) is the system
of ODEs of the form
\begin{equation} \label {e2}
\begin{array}{l}
   \rho _1 \,Z'' - k e^{ -\frac{d}{{2\rho _2 }}\zeta } \left(\beta \cos (\frac{{\mu \zeta }}{{2\rho _2 }}) +\gamma \sin (\frac{{\mu \zeta }}{{2\rho _2 }})\right) - \alpha d = 0, \\
    \rho _2 \,W'' + d\,W' + k\,W = 0. \\ \nonumber
\end{array}
\end{equation}
Solving this system yields the following solution
\begin{equation} \label {34}
\begin{array}{ll}
  \varphi (t,x) =& \left( \left( {\beta\,\mu ^2- 2\,\beta \,k\rho _2  - d\gamma \,\mu } \right)\cos \left( {\,\frac{{\mu \,t}}{{2\rho _2 }}} \right)+\left( {\gamma \,\mu ^2  - 2\,\gamma \,k\rho _2 }\right)\sin \left( {\,\frac{{\mu \,t}}{{2\rho _2 }}} \right) \right){\rm{e}}^{ - \frac{{dt}}{{2\rho _2 }}} \\
  & - \frac{{\beta\mu d}}{{2\rho _2 }}\,\,\left( { - \alpha \,kx^2 \rho _1 \rho _2+ d} \right)t\sin \left( {\,\frac{{\mu \,t}}{{2\rho _2 }}} \right)- \left( {\alpha \,kdt^2  + 2\,c_1 \,\rho _1 kt + 2\,c_2 } \right) ,\\
  \psi (t,x) =& \left( {\left( {\beta x + c_3 } \right)\cos \left({\,\frac{{\mu \,t}}{{2\rho _2 }}} \right) + \left( {\gamma x + c_4} \right)\sin \left( {\,\frac{{\mu \,t}}{{2\rho _2 }}} \right)}\right){\rm{e}}^{ - \,\frac{{dt}}{{2\rho _2 }}}  + \alpha\,x\left( {\frac{d}{k} - t} \right),\\
\end{array}
\end{equation}
where $d=\sqrt{k\rho_2-\mu^2}$.
\end{example}
   \begin{center}
\begin{table}[H]  \label{tab:7}
\centering
\begin{scriptsize}
 \caption {Reduction using one-dimensional optimal system (\ref{28}) with $d=2\sqrt{k\rho_2}$}
\begin{tabular}{ c c c }
  \hline
 \rowcolor[gray] {0.9}  Generators in (\ref{30})   & Invariant & The reduced system  \\
\rowcolor[gray] {0.9} &variables&\\ \hline
$X^1=X_1+\alpha X_2 +\beta X_6, $ &A & $ (k - \alpha ^2 \rho _1)\,Z'' + k\,W' - \beta \rho _1 \,\zeta \, = \,0, $\\
        $\alpha \in \mathbb{R},\beta \neq 0.$ & & $(\chi '(W') - \alpha ^2 \rho _2 )\,W'' - k\,Z' + \alpha d\,W' -k\,W - 3\beta \rho _2 \, = \,0. $ \\
   \hline
$X^2=X_1+\alpha X_2 +\beta X_4$& B & $(k - \alpha ^2 \rho _1)\,Z'' + k\,W' - \beta \rho _1 \, = \,0,$ \\
         $\alpha,\beta \in \mathbb{R}.$   & &$(\chi '(W') - \alpha ^2 \rho _2 )\,W'' - kZ' + \alpha d\,W' - k\,W= 0.$  \\
   \hline
$X^3=X_2 +\alpha X_6+\beta X_7 +\gamma X_8$ & C &$ \rho _1 \,Z'' -k(\gamma \zeta  + \beta )e^{ -  \frac{d}{{2 \rho_2 }} \zeta} {\rm{ }} - \alpha d\, = \,0,$  \\
         $\alpha \neq 0,\beta,\gamma \in \mathbb{R}.$   &  &$\rho _2 \,W'' + d\,W' + k\,W\, = \,0. $  \\
   \hline
$X^4=X_2 +\alpha X_5+\beta X_8$  & D & $ \rho _1 \,Z'' - \beta k\,\zeta \,e^{ - \frac{d}{{2 \rho_2 }} \zeta } \, = \,0,$ \\
         $\alpha \in \mathbb{R},\beta \neq 0.$    &  & $\rho _2 \,W'' + d\,W' + k\,W\, = \,0. $ \\
   \hline
$X^5=\alpha X_2 +\beta X_5+ X_7 $  & E & $  \rho _1 \,Z'' -\frac{k}{\alpha }\,e^{ -\frac{d}{{2 \rho_2 }} \zeta } \, =\,0, $ \\
         $\alpha \neq 0,\beta  \in \mathbb{R}.$   &   & $ \rho _2 \,W'' + d\,W' + k\,W\, = \,0.$ \\
   \hline
$X^6=X_2 +\alpha X_5$  & F & $Z''\, = \,0,$ \\
           $\alpha \in \mathbb{R}.$    & & $\rho _2 \,W'' + d\,W' + k\,W\, = \,0.$ \\
   \hline
\end{tabular}
\end{scriptsize}
\end{table}
\begin{scriptsize}
\begin{equation}\label{}
\begin{array}{ll|l|l}
     A:& \phi (t,x)\, = \,Z(\zeta ) + \frac{\beta }{2}t^2 (x - \frac{\alpha}{3}t), & \psi (t,x)\, = \,W(\zeta ) - \frac{\beta }{2}t^2  + \frac{{\beta d}}{k}t, & \zeta \, = \,x - \alpha t. \\
     B:&\phi (t,x)\, = \,Z(\zeta ) + \frac{\beta }{2}t^2 ,  & \psi (t,x)\, = \,W(\zeta ), &\zeta \, = \,x - \alpha t. \\
     C:& \phi (t,x)\, = \,Z(\zeta ) + \frac{\alpha }{2}tx^2 , &\psi (t,x)\, = \,W(\zeta ) + \alpha x(\frac{d}{k} - t) + x(\beta+ \gamma t)e^{ - \frac{d}{{2 \rho_2 }} t} , & \zeta \, = \,t. \\
     D:& \phi (t,x)\, = \,Z(\zeta ) + \frac{\alpha }{2}x^2 , & \psi (t,x) = W(\zeta ) - x(\alpha  - \beta t\,e^{ - \frac{d}{{2 \rho_2 }} t} ), & \zeta \, = \,t. \\
     E:& \phi (t,x)\, = \,Z(\zeta ) + \frac{\beta }{{2\alpha }}x^2 , & \psi (t,x)\, = \,W(\zeta ) - \frac{x}{\alpha }(\beta  - \,e^{ -\frac{d}{{2 \rho_2 }} t} ), & \zeta \, = \,t.  \\
     F:& \phi (t,x)\, = \,Z(\zeta ) + \frac{\alpha }{2}x^2 , & \psi (t,x)\, = \,W(\zeta ) - \alpha x, & \zeta \, = \,t. \\
\nonumber
    \end{array}
\end{equation}
\end{scriptsize}
\end{center}
\begin{center}
\begin{table}[H]  \label{tab:8}
\centering
\begin{scriptsize}
 \caption {Reduction using one-dimensional optimal system (\ref{29}) with $d=\sqrt{4k\rho_2+\lambda^2}$}
\begin{tabular}{ c c c }
  \hline
 \rowcolor[gray] {0.9}  Generators in (\ref{30})   & Invariant & The reduced system  \\
\rowcolor[gray] {0.9} &variables&\\ \hline
  $X^1=X_1+\alpha X_2 +\beta X_6, $  &A& $(k - \alpha ^2 \rho _1 )\,Z'' + k\,W' - \beta \rho _1 \,\zeta \, =\,0, $\\
        $\alpha \in \mathbb{R},\beta \neq 0.$   & & $(\chi '(W') - \alpha ^2 \rho _2 )\,W'' - k\,Z' + \alpha d\,W' -k\,W - \beta (3\rho _2  - \frac{{\lambda ^2 }}{k})\, = \,0. $  \\
    \hline
        $X^2=X_1+\alpha X_2 +\beta X_4, $ &B & $(k - \alpha ^2 \rho _1 )\,Z'' + k\,W' - \beta \rho _1 \, = \,0,$\\
        $\alpha,\beta \in \mathbb{R}.$  & & $ (\chi '(W') - \alpha ^2 \rho _2 )\,W'' - k\,Z' + \alpha d\,W' -k\,W = 0.$ \\
    \hline
       $X^3=X_2+\alpha X_6 +\beta X_7+\gamma X_8 $  &C& $\rho _1 \,Z'' - k\,e^{ -\frac{d}{{2\rho _2 }}\zeta }(\beta \cosh (\frac{{\lambda \zeta }}{{2\rho _2}}) + \gamma \sinh (\frac{{\lambda \zeta }}{{2\rho _2 }}))  - \alpha d\, = \,0, $\\
          $\alpha \neq 0,\beta,\gamma \in \mathbb{R}.$& & $ \rho _2 \,W'' + d\,W' + k\,W\, = \,0.$  \\
    \hline
       $X^4=X_2+\alpha X_5 +\beta X_7 $ &D & $\rho _1 Z'' - \beta k\,e^{ - \frac{d}{{2\rho _2 }}\zeta } \,\cosh(\frac{{\lambda \zeta }}{{2\rho _2 }})\, = \,0, $\\
        $\alpha \in \mathbb{R},\beta \neq 0.$   & & $\rho _2 \,W'' + d\,W' + k\,W\, = \,0. $  \\
    \hline
     $X^5=X_2+\alpha X_5 +\beta X_8 $ &E & $ \rho _1 \,Z'' - \beta k\,e^{ - \frac{d}{{2\rho _2 }}\zeta } \sinh(\frac{{\lambda \zeta}}{{2\rho _2 }})\, = \,0,$\\
        $\alpha \in \mathbb{R},\beta \neq 0.$  & & $\rho _2 \,W'' + d\,W' + k\,W\, = \,0. $  \\
    \hline
     $X^6=\alpha X_2+\beta X_5 + X_7+ X_8 $ &F & $\rho _1 \,Z'' - \frac{k}{\alpha }\,e^{ -\frac{d}{{2\rho _2 }}\zeta }(\cosh (\frac{{\lambda \zeta}}{{2\rho_2 }}) + \sinh (\frac{{\lambda \zeta}}{{2\rho _2 }})) \, = \,0, $\\
        $\alpha \neq 0,\beta \in \mathbb{R}.$  & & $ \rho _2 \,W'' + d\,W' + k\,W\, = \,0.$  \\
    \hline
     $X^7=\alpha X_2+\beta X_5 + X_7- X_8 $ &G & $\rho _1 \,Z'' - \frac{k}{\alpha }e^{ -\frac{d}{{2\rho _2 }}\zeta }(\cosh (\frac{{\lambda \zeta}}{{2\rho_2 }}) - \sinh (\frac{{\lambda \zeta}}{{2\rho _2 }})) \, = \,0, $\\
        $\alpha \neq 0,\beta \in \mathbb{R}.$ &  & $\rho _2 \,W'' + d\,W' + k\,W\, = \,0. $ \\
    \hline
     $X^8=X_2+\alpha X_5 $ &H & $Z''\, = \,0,$\\
           $\alpha \in \mathbb{R}.$& & $\rho _2 \,W'' + d\,W' + k\,W\, = \,0. $ \\
    \hline
\end{tabular}
\end{scriptsize}
\end{table}
\begin{scriptsize}
\begin{equation}\label{}
\begin{array}{ll|l|l}
     A: &\phi (t,x)\, = \,Z(\zeta ) - \frac{\beta }{6}t^2 (\alpha t -3x),\,&\psi (t,x)\, = \,W(\zeta ) - \frac{\beta }{2}t(t -\frac{{2d}}{k}),\,&\zeta \, = \,x - \alpha t. \\
     B: & \phi (t,x)\, = \,Z(\zeta ) + \frac{\beta }{2}t^2 ,\,&\psi (t,x)\, =\,W(\zeta ),\,&\zeta \, = \,x - \alpha t.\\
     C: &\phi (t,x)\, = \,Z(\zeta ) + \frac{\alpha }{2}tx^2 ,\,&\psi (t,x)\,= \,W(\zeta ) + \alpha x(\frac{d}{k} - t) + xe^{ - \frac{d}{{2\rho_2 }}t} (\beta \cosh (\frac{{\lambda t}}{{2\rho _2 }}) + \gamma\sinh (\frac{{\lambda t}}{{2\rho _2 }})),\, & \zeta \, = \,t.\\
     D: &\phi (t,x)\, = \,Z(\zeta ) + \frac{\alpha }{2}x^2 ,\,&\psi (t,x)\,= \,W(\zeta ) - x(\,\alpha  - \beta e^{ - \frac{d}{{2\rho _2 }}t}\cosh (\frac{{\lambda t}}{{2\rho _2 }})\,),\,&\zeta \, = \,t. \\
     E: &\phi (t,x)\, = \,Z(\zeta ) + \frac{\alpha }{2}x^2 ,\,&\psi (t,x)\,= \,W(\zeta ) - x(\,\alpha  - \beta e^{ - \frac{d}{{2\rho _2 }}t}\sinh (\frac{{\lambda t}}{{2\rho _2 }})\,),\,&\zeta \, = \,t.  \\
     F: &\phi (t,x)\, = \,Z(\zeta ) + \frac{\beta }{{2\alpha }}x^2 ,\, & \psi(t,x)\, = \,W(\zeta ) + \frac{x}{\alpha }e^{ - \frac{d}{{2\rho _2}}t} (\cosh (\frac{{\lambda t}}{{2\rho _2 }}) + \sinh(\frac{{\lambda t}}{{2\rho _2 }}) - \beta ),\,&\zeta \, = \,t. \\
     G: &\phi (t,x)\, = \,Z(\zeta ) + \frac{\beta }{{2\alpha }}x^2 ,\, & \psi(t,x)\, = \,W(\zeta ) + \frac{x}{\alpha }e^{ - \frac{d}{{2\rho _2}}t} (\cosh (\frac{{\lambda t}}{{2\rho _2 }}) - \sinh(\frac{{\lambda t}}{{2\rho _2 }}) - \beta ),\, & \zeta \, = \,t.\\
     H: &\phi (t,x)\, = \,Z(\zeta ) + \frac{\alpha }{2}x^2 ,\,&\psi (t,x) =W(\zeta ) - \alpha x,\,&\zeta \, = \,t.  \\ \nonumber
    \end{array}
\end{equation}
\end{scriptsize}
\end{center}
\begin{center}
\begin{table}[H]  \label{tab:9}
\centering
\begin{scriptsize}
 \caption {Reduction using one-dimensional optimal system (\ref{30}) with $d=\sqrt{4k\rho_2-\mu^2}$}
\begin{tabular}{ c c c c }
  \hline
\rowcolor[gray] {0.9}  Generators in (\ref{30})   & Invariant & The reduced system  \\
\rowcolor[gray] {0.9} &variables&\\ \hline
  $X^1=X_1+\alpha X_2 +\beta X_6 $  &A& $(k - \alpha ^2 \rho _1 )\,Z'' + k\,W' - \beta \rho _1 \,\zeta \, =\,0,$\\
         $\alpha \in \mathbb{R},\beta \neq 0.$ & & $(\chi '(W') - \alpha ^2 \rho _2 )\,W'' - k\,Z' + \alpha d\,W' -k\,W + \beta (\frac{{\mu ^2 }}{k}{\rm{ - 3}}\rho _2 )\, = \,0.$  \\
    \hline
       $X^2=X_1+\alpha X_2 +\beta X_4$ &B &$(k - \alpha ^2 \rho _1 )\,Z'' + k\,W' - \beta \rho _1 \, = \,0,$ \\
          $\alpha,\beta \in \mathbb{R}.$   & &$(\chi '(W') - \alpha ^2 \rho _2 )\,W'' - k\,Z' + \alpha d\,W' -k\,W = 0.$  \\
         \hline
       $X^3=X_2 +\alpha X_6+\beta X_7 +\gamma X_8$& C &$\rho _1 \,Z'' - k e^{ -\frac{d}{{2\rho _2 }}\zeta } \left(\beta \cos (\frac{{\mu \zeta }}{{2\rho _2 }}) +\gamma \sin (\frac{{\mu \zeta }}{{2\rho _2 }})\right) - \alpha d = 0,$  \\
           $\alpha \neq 0,\,\beta,\gamma \in \mathbb{R}.$  & &$\rho _2 \,W'' + d\,W' + k\,W = 0.$  \\
          \hline
       $X^4=X_2 +\alpha X_5+\beta X_7 +\gamma X_8$  &D & $\rho _1 \,Z'' - k e^{ -\frac{d}{{2\rho _2 }}\zeta }\left(\beta \cos (\frac{{\mu \zeta }}{{2\rho _2 }}) +\gamma \sin (\frac{{\mu \zeta }}{{2\rho _2 }})\right) \, = \,0,$ \\
           $\alpha,\beta,\gamma \in \mathbb{R}.$  &  & $\rho _2 \,W'' + d\,W' + k\,W = 0.$ \\
            \hline
\end{tabular}
\end{scriptsize}
\end{table}
\begin{scriptsize}
\begin{equation}\label{}
\begin{array}{ll|l|l}
     A:& \phi(t,x) \, = \,Z(\zeta ) - \frac{\beta }{6}t^2 (\alpha t - 3x),&\, \psi(t,x) \, = \,W(\zeta ) - \frac{\beta }{2}t(t - \frac{{2d}}{k}),& \zeta \, = \,x - \alpha t. \\
     B:& \phi(t,x) \, = \,Z(\zeta ) + \frac{\beta }{2}t^2 ,& \psi(t,x) \, = \,W(\zeta ),& \zeta \, = \,x - \alpha t.\\
     C:& \phi (t,x)\, = \,Z(\zeta ) + \frac{\alpha }{2}tx^2 ,& \psi(t,x) \, = \,W(\zeta ) + \alpha x(\frac{d}{k} - t) + xe^{- \frac{d}{{2\rho _2 }}t}\left(\beta \cos (\frac{{\mu t}}{{2\rho _2 }}) + \gamma \sin (\frac{{\mu t}}{{2\rho _2 }})\right) ,& \zeta \, = \,t.{\rm{ }}\\
     D:& \phi (t,x)\, = \,Z(\zeta ) + \frac{\alpha }{2}x^2 ,& \psi(t,x) \, = \,W(\zeta ) - x\,\alpha  + xe^{ - \frac{d}{{2\rho _2 }}t}\left(\beta \cos (\frac{{\mu t}}{{2\rho _2 }}) - \gamma \sin (\frac{{\mu t}}{{2\rho _2}})\,\right) ,& \zeta \, = \,t. \\  \nonumber
    \end{array}
\end{equation}
\end{scriptsize}
\end{center}

\section{Discussion and concluding remarks}
The complete Lie symmetry classification of a non-linear
Timoshenko system of PDEs with frictional damping term in
rotational angle is performed. The classification is related to
the arbitrary dependence on the rotation moment $\chi(\psi_x)$. A
Lie symmetry analysis is performed in three cases for non-linear
rotational moment. The three cases depend on the sign of the
parameter $d^2-4 k\rho_2$. The one-dimensional optimal system is
derived for each one of the three cases. All possible invariant
forms and their corresponding reductions for each vector field in
the optimal systems are found. These reductions to systems of ODEs
are given in Table 4, Table 5 and Table 6. They are described by
optimal reduction where all non-similar invariant solutions under
symmetry transformations can be given from the solution of these
reduced system of ODEs.
\subsection*{Acknowledgments}
The authors are grateful to King Fahd University of Petroleum \&
Minerals, Saudi Arabia for supporting and providing research
facilities.
\section*{Appendix: The Tree diagrams }
\begin{center}
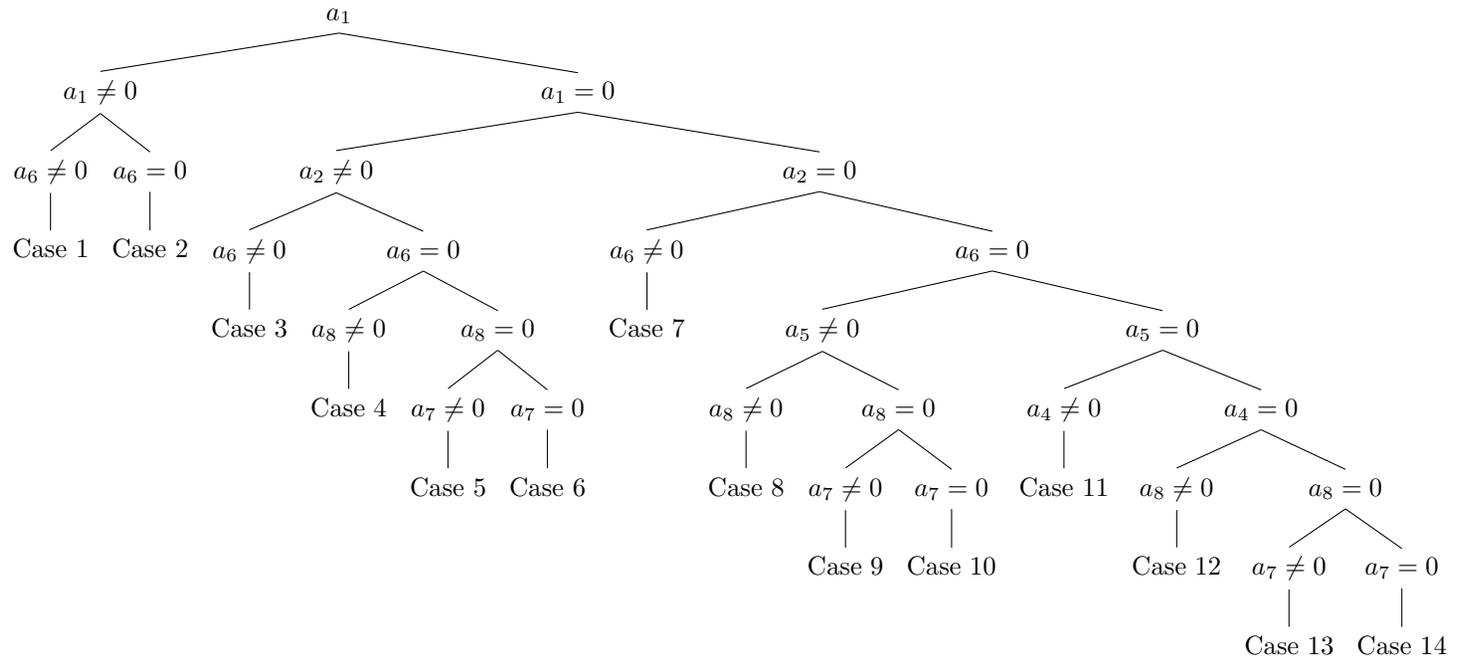
\begin{sidewaysfigure}

\begin{figure}[H]
\caption{Tree diagram 1}
\begin{tikzpicture}[scale=1]

\Tree [.{$a_1$}    [.{$a_1\neq0$} [ .{$a_6\neq0$} Case~1 ] [
.{$a_6 =0$} Case~2 ] ]
   [ .{$a_1 = 0$} [ .{$a_2\neq 0$}  [ .{$a_6\neq 0$} Case~3 ] [ .{$a_6 = 0$} [ .{$a_8 \neq 0$} Case~4 ]
   [ .{$a_8 = 0$} [ .{$a_7 \neq 0$} Case~5 ][ .{$a_7 = 0$} Case~6 ] ] ]  ] [ .{$a_2 = 0$} [ .{$a_6\neq0$} Case~7 ]
   [ .{$a_6=0$} [ .{$a_5 \neq 0$} [ .{$a_8 \neq 0$} Case~8 ][ .{$a_8 = 0$} [ .{$a_7 \neq 0$} Case~9 ]
   [ .{$a_7 = 0$} Case~10 ] ] ][ .{$a_5=0$} [ .{$a_4 \neq 0$} Case~11 ][ .{$a_4=0$} [ .{$a_8\neq0$} Case~12 ]
   [ .{$a_8=0$} [ .{$a_7 \neq 0$} Case~13 ][ .{$a_7=0$} Case~14 ] ] ] ] ] ] ] ]

\end{tikzpicture}
\end{figure}

\end{sidewaysfigure}
\end{center}
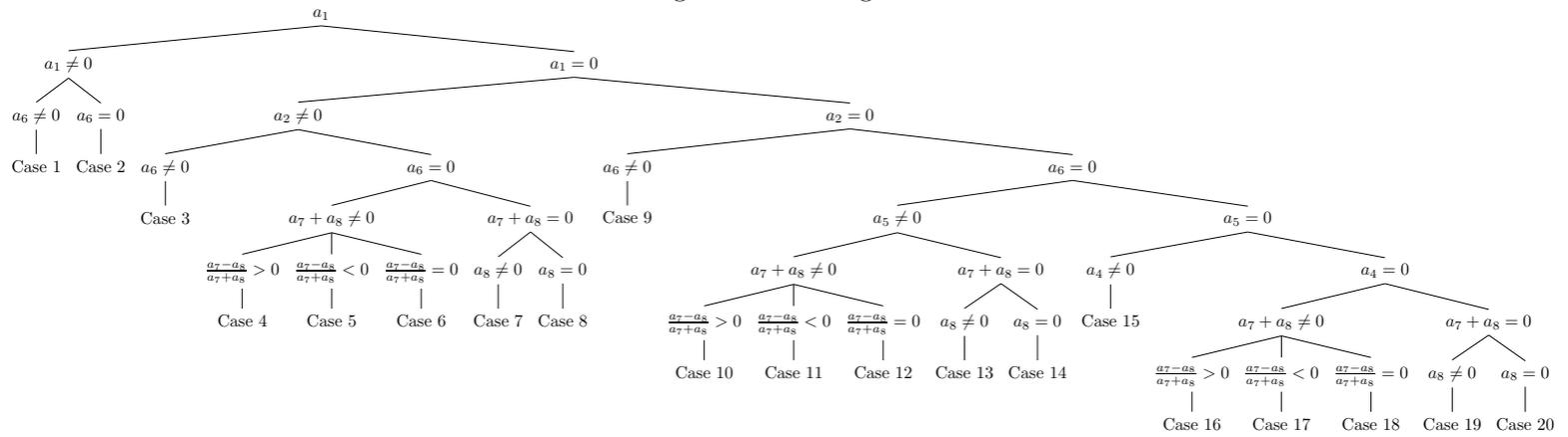
\begin{sidewaysfigure}
\begin{center}
\begin{figure}[H]
\caption{Tree diagram 2}
\begin{tikzpicture}[scale=.65]

\Tree [.{$a_1$}    [.{$a_1\neq 0$} [ .{$a_6\neq0$} Case~1 ] [
.{$a_6 = 0$} Case~2 ] ]
   [ .{$a_1 = 0$} [ .{$a_2\neq 0$}  [ .{$a_6\neq 0$} Case~3 ] [ .{$a_6 = 0$} [ .{$a_7+a_8 \neq 0$}
   [ .{$ \frac{a_7-a_8}{a_7+a_8}>0$} Case~4 ] [ .{$\frac{a_7-a_8}{a_7+a_8}<0$} Case~5 ] [ .{$\frac{a_7-a_8}{a_7+a_8}=0$} Case~6 ]  ]
   [ .{$a_7+a_8 =0$} [ .{$a_8 \neq 0$} Case~7 ][ .{$a_8 = 0$} Case~8 ] ] ]  ]
    [ .{$a_2 = 0$} [ .{$a_6\neq0$} Case~9 ]
   [ .{$a_6=0$} [ .{$a_5 \neq 0$} [ .{$a_7+a_8 \neq 0$}  [ .{$\frac{a_7-a_8}{a_7+a_8}>0$} Case~10 ]
   [ .{$\frac{a_7-a_8}{a_7+a_8}<0$} Case~11 ][ .{$\frac{a_7-a_8}{a_7+a_8}=0$} Case~12 ]] [ .{$a_7+a_8= 0$}
    [ .{$a_8 \neq 0$} Case~13 ]
   [ .{$a_8 = 0$} Case~14 ] ] ][ .{$a_5=0$} [ .{$a_4 \neq 0$} Case~15 ]  [ .{$a_4=0$} [ .{$a_7+a_8 \neq 0$}
    [ .{$\frac{a_7-a_8}{a_7+a_8}>0$} Case~16 ][ .{$\frac{a_7-a_8}{a_7+a_8}<0$} Case~17 ][ .{$\frac{a_7-a_8}{a_7+a_8}=0$} Case~18 ] ]
   [ .{$a_7+a_8= 0$} [ .{$a_8 \neq 0$} Case~19 ][ .{$a_8=0$} Case~20 ] ] ] ] ] ] ]
   ]]
\end{tikzpicture}
\end{figure}
\end{center}
\end{sidewaysfigure}
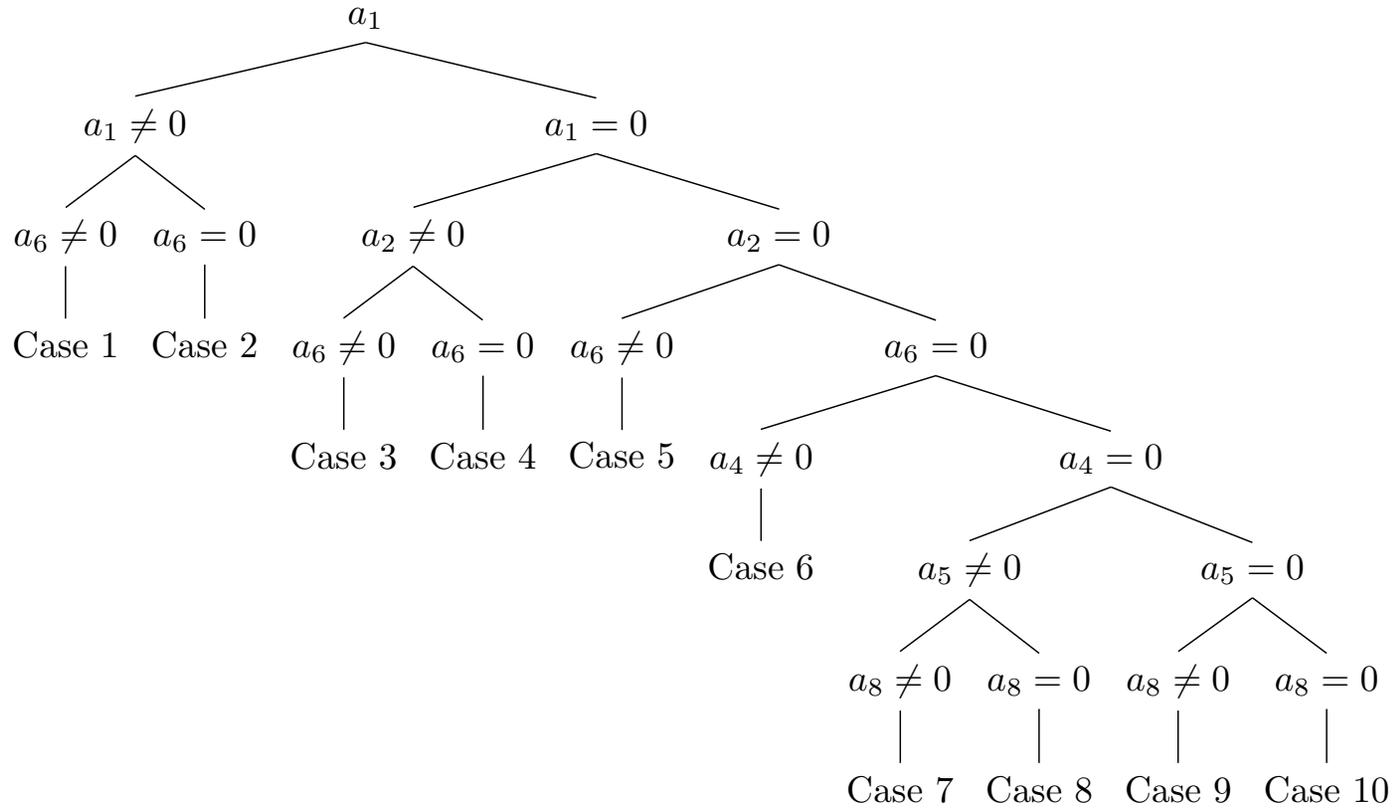
\begin{sidewaysfigure}
\begin{center}
\begin{figure}[H]
\caption{Tree diagram 3}
\begin{tikzpicture}[scale=1.4]
\Tree [.{$a_1$}    [.{$a_1\neq0$} [ .{$a_6\neq0$} Case~1 ] [
.{$a_6=0$} Case~2 ] ][ .{$a_1 = 0$} [ .{$a_2\neq 0$}  [ .{$a_6\neq
0$} Case~3 ] [ .{$a_6 = 0$} Case~4 ]  ] [ .{$a_2 = 0$} [
.{$a_6\neq0$} Case~5 ] [ .{$a_6=0$} [ .{$a_4 \neq 0$} Case~6 ] [
.{$a_4= 0$} [ .{$a_5 \neq 0$} [ .{$a_8 \neq 0$} Case~7 ][
.{$a_8=0$} Case~8 ] ][ .{$a_5 = 0$} [ .{$a_8 \neq 0$} Case~9 ][
.{$a_8 = 0$} Case~10 ] ] ] ] ] ] ]
\end{tikzpicture}
\end{figure}
\end{center}
\end{sidewaysfigure}

\end{document}